\documentstyle[12pt]{article}
\begin{document}
\newtheorem{Def}{Definition}[section]
\newtheorem{thm}{Theorem}[section]
\newtheorem{lem}{Lemma}[section]
\newtheorem{rem}{Remark}[section]
\newtheorem{question}{Question}[section]
\newtheorem{prop}{Proposition}[section]
\newtheorem{cor}{Corollary}[section]
\newtheorem{clm}{Claim}[section]
\newtheorem{step}{Step}[section]
\newtheorem{sbsn}{Subsection}[section]
\newtheorem{conj}{Conjecture}[section]
\title
{On some conformally invariant fully
nonlinear equations, Part II: Liouville, Harnack and Yamabe }
\author{\ Aobing Li
\ \ \ \ \& \ \ \ YanYan Li
\\ Department of Mathematics\\ Rutgers University\\
110 Frelinghuysen Rd.\\
Piscataway, NJ 08854
}
\date{}
\maketitle
\input { amssym.def}
\begin{abstract}
The Yamabe problem concerns finding a conformal metric on a given
closed Riemannian manifold so that it has constant scalar
curvature. This paper concerns mainly a fully nonlinear version of the Yamabe
problem and the corresponding Liouville type problem.
\end{abstract}
\setcounter{section}{0}
\section{Introduction}

\bigskip

Let $(M,g)$ be an
 $n$-dimensional compact smooth Riemannian manifold (without
boundary).
For $n=2$, we know from the uniformization theorem of
Poincar\'e that there exist metrics that are pointwise
conformal to $g$ and have constant Gauss curvature.
For $n\ge 3$, the well-known  Yamabe
conjecture
states that
there exist metrics which are pointwise conformal to $g$ and
have constant
scalar curvature.  The  Yamabe conjecture is
proved through the work of Yamabe \cite{Y}, Trudinger
\cite{Tr},
Aubin \cite{A} and  Schoen \cite{S0}.
The Yamabe and related problems
have attracted much attention in the last
30 years or so, see, e.g.  \cite{SY6}, \cite{Aubin},
and the references therein.
Important methods and techniques in overcoming
loss of compactness  have been
developed in such studies which also play
important roles in the research of
other areas of mathematics.
For $n\ge 3$, let $\hat g=u^{ \frac 4{n-2} }g$ where
$u$ is  some
positive function  on $M$. The scalar curvature $R_{\hat g}$ of
$\hat g$ can be calculated as
$$
R_{\hat g}=u^{ -\frac {n+2}{n-2} }
\bigg( R_g u-\frac {4(n-1) }{ n-2} \Delta_g u\bigg),
$$
where $R_g$ and $\Delta_g$ denote
respectively the scalar curvature and the Laplace-Beltrami
operator  of $g$.
The Yamabe conjecture is therefore equivalent to the existence
of a positive solution of
\[
-L_g u=\bar R u^{\frac{n+2}{n-2}}\quad\mbox{on}~M,
\]
where  $L_g:=\Delta_g-\frac{n-2}{4(n-1)}R_g$ is the
conformal Laplacian of $g$, and $\bar R=0$ or
$\pm 2(n-1)$.
The Yamabe problem can be divided into three cases --- positive case, zero case
and negative case ---
according to the signs of the first eigenvalue
of $-L_g$. Making a conformal change of metrics
$\tilde g=\varphi^{ \frac 4{n-2} }g$, where $\varphi$
is
 a positive eigenfunction of $-L_g$ associated with the first eigenvalue,
we are led to the following three cases: $R_g>0$ on $M$,
$R_g\equiv 0$ on $M$ and $R_g<0$ on $M$.
The positive case, i.e. $R_g>0$, is much more difficult.

Let
$$
A_g:=\frac{1}{n-2}(Ric_g-\frac{R_g}{2(n-1)}g)
$$
denote the Schouten tensor of $g$, where $Ric_g$ denotes the Ricci tensor
of $g$.   We use $\lambda(A_g)=(\lambda_1(A_g),
\cdots, \lambda_n(A_g))$ to denote the eigenvalues of
$A_g$ with respect to $g$. Clearly
$$
\sum_{i=1}^n \lambda_i(A_g)=\frac 1{2(n-1)}R_g.
$$

Let
$$
V_1=\{\lambda\in \Bbb R^n\ | \
\sum_{i=1}^n \lambda_i>1\},
$$
and let
$$
\Gamma (V_1)=\{s\lambda|~s>0,~\lambda\in V_1\}
$$
be the cone with vertex at the origin generated by $V_1$.

The Yamabe problem in the positive case can be reformulated as
follows:  Assuming $\lambda(A_g)\in \Gamma(V_1)$, then there
exists a Riemannian metric $\hat g$ which is pointwise conformal
to $g$ and satisfies $\lambda(A_{\hat g})\in \partial V_1$ on $M$.

In general, let $V$ be an  open convex 
subset of $\Bbb R^n$ which is symmetric with respect
to the coordinates, i.e., $(\lambda_1, \cdots, \lambda_n)
\in V$ implies $(\lambda_{i_1}, \cdots, \lambda_{i_n})
\in V$ for any permutation $(i_1, \cdots, i_n)$
of $(1,\cdots, n)$.  We assume that $\emptyset\neq\partial V$ is
in $C^{2,\alpha}$ for some $\alpha\in (0,1)$ in the sense that $\partial V$ can be
represented as the graph of some $C^{2,\alpha}$ function near every point.
For $\lambda\in \partial V$, let $\nu(\lambda)$
denote the inner unit normal of $\partial V$.  We further assume
that
\begin{equation}
\nu(\lambda)\in \Gamma_n:= \{\lambda\in R^n|\lambda_i>0,\forall 1\le
i\le n\},\qquad\forall\ \lambda\in \partial V,
\label{nu1}
\end{equation}
and
\begin{equation}
\nu(\lambda)\cdot \lambda>0, \qquad \forall\ \lambda\in \partial V.
\label{nu2}
\end{equation}

Let 
\begin{equation}
\Gamma(V):=\{s\lambda\ |\ \lambda \in V, \ 0<s<\infty\}.
\label{cone}
\end{equation}
be the (open convex) cone with vertex at the origin generated
by $V$.

Our first theorem establishes the existence and compactness of solutions
to a fully nonlinear version of the Yamabe problem
on locally
conformally flat manifolds.
A Riemannian manifold $(M^n, g)$ is called locally conformally flat
if near every point of $M$ the metric can be represented in
some local coordinates as $g=e^{\psi(x)}\sum_{i=1}^n
(dx^i)^2$.

\begin{thm} For $n\ge 3, \alpha\in
(0,1)$, we assume  that
 $V$  is a symmetric open convex subset of $R^n$, with
$\emptyset\neq \partial V\in C^{4,\alpha}$,
satisfying
(\ref{nu1}) and (\ref{nu2}). Let  $(M^n,g)$ be a compact,
smooth, connected,
locally conformally flat Riemannian manifold of dimension
$n$ satisfying 
$$
\lambda(A_g)\in \Gamma(V),\qquad \mbox{on}\ M^n.
$$
Then there exists a positive function $u\in C^{4,\alpha}(M^n)$ such that
the conformal metric $\hat g=u^{ \frac 4{n-2} }g$
satisfies
\begin{equation}
\lambda(A_{\hat g})\in \partial V, \qquad \mbox{on}\ M^n.
\label{yamabe1}
\end{equation}
Moreover, if $(M^n,g)$ is not conformally diffeomorphic
to the standard $n-$sphere, then
all positive solutions of (\ref{yamabe1}) satisfy
\[
 \|u\|_{ C^{4,\alpha}(M^n,g) }+  \|\frac 1u\|_{ C^{4,\alpha}(M^n,g) }
\le C, \qquad \mbox{on}\ M^n,
\]
where $C$ is some positive constant depending only on
$(M^n,g)$, $V$ and $\alpha$.
\label{yamabe}
\end{thm}

\begin{rem}  Presumably, the existence of a $C^{2,\alpha}$ 
solution of (\ref{yamabe1}) should hold under the weaker smoothness
hypothesis $\partial V\in C^{2,\alpha}$.  We prove this
under an additional hypothesis that $V$ is strictly convex, i.e.,
principal curvatures of $\partial V$ are positive everywhere. See
Appendix B.
\label{remark1-1}
\end{rem}

We propose the following
\begin{conj}
Assume that $V$ is an open symmetric convex subset of
$\Bbb {R}^n$, with $\emptyset\neq\partial V\in C^\infty$, satisfying
(\ref{nu1}) and (\ref{nu2}).  Let
$(M^n,g)$ be
 a compact smooth Riemannian manifold of dimension $n\ge 3$
 satisfying
\[
\lambda(A_g)\in \Gamma(V),\qquad \mbox{on}\ M^n.
\]
Then there exists a smooth positive function
$u\in C^\infty(M^n)$ such that the conformal metric
 $\hat g=u^{ \frac 4{n-2} }g$ satisfies
\begin{equation}
\lambda(A_{\hat g})\in \partial V, \qquad \mbox{on}\ M^n.
\label{Q2}
\end{equation}
\label{question1}
\end{conj}

For $V=V_1$, it is the Yamabe problem in the positive case.
In general, the equation of $u$ is a fully nonlinear
elliptic equation of second order, and therefore
the problem can be viewed as a fully nonlinear version of
the Yamabe problem.

The fully nonlinear
version of the Yamabe problem has the following equivalent formulation.
The equivalence of the two formulations is shown in Appendix B.
 
Let
\begin{equation}
\Gamma\subset \Bbb R^n \ \mbox{be
an open convex symmetric cone  with vertex at the origin }
\label{gamma0}
\end{equation}
satisfying
\begin{equation}
\Gamma_n
\subset \Gamma\subset \Gamma_1:=\{\lambda\in\Bbb{R}^n|\sum\limits_i
\lambda_i>0\}.  
\label{gamma1}
\end{equation}
Naturally, $\Gamma$ being symmetric  
means $(\lambda_1, \lambda_2,\cdots, \lambda_n)\in
\Gamma$ implies $(\lambda_{i_1}, \lambda_{i_2}, \cdots , \lambda_{i_n})
\in \Gamma$ for any permutation 
 $(i_1, i_2, \cdots , i_n)$ of
$(1,2,\cdots, n)$.

For $\alpha\in (0,1)$, let
\begin{equation}
f\in C^{4,\alpha}(\Gamma)\cap
C^0(\overline \Gamma)\ 
\mbox{be  concave and symmetric in}\ \lambda_i,
\label{hypo0}
\end{equation}
satisfying
\begin{equation}
f|_{\partial\Gamma}=0,\qquad
\nabla f\in\Gamma_n\ \mbox{on}\ \Gamma,
\label{aa5}
\end{equation}
and
\begin{equation}
\lim_{s\to\infty}f(s\lambda)=\infty,\qquad
\forall\ \lambda\in \Gamma.
\label{aa6}
\end{equation}

Conjecture \ref{question1} is equivalent to

\noindent{\bf Conjecture $\ref{question1}^{\prime}$.}\
{\it  Assume that  $(f,\Gamma)$ satisfies
(\ref{gamma0}), (\ref{gamma1}),
 (\ref{hypo0}), (\ref{aa5})
and (\ref{aa6}). Let $(M^n, g)$ be
 a compact smooth Riemannian manifold of dimension $n\ge 3$,
 satisfying
  $\lambda
(A_g)\in\Gamma$ on $M^n$. Then there exists  a smooth positive function
$u\in C^\infty(M^n)$ such that the conformal metric
 $\hat g=u^{ \frac 4{n-2} }g$ satisfies
\begin{equation}
f\left(\lambda(A_{\hat g})\right)=1,
\quad \lambda(A_{\hat g})\in \Gamma, \qquad \mbox{on}\ M^n.
\label{ghat2}
\end{equation}
}

 Theorem~\ref{yamabe} is equivalent to

\noindent{\bf Theorem \ref{yamabe}$'$}.\
{\it  For $n\ge 3, \alpha\in
(0,1)$, we assume  that  $(f,\Gamma)$ satisfies
(\ref{gamma0}), (\ref{gamma1}), 
 (\ref{hypo0}), (\ref{aa5})
and (\ref{aa6}).  Let  $(M^n,g)$ be a compact,
smooth, connected,
locally conformally flat Riemannian manifold of dimension
$n$ satisfying
  $\lambda
 (A_g)\in\Gamma$ on $M^n$.
Then  there exists a positive function $u\in C^{4,\alpha}(M^n)$ such that
the conformal metric $\hat g=u^{ \frac 4{n-2} }$
satisfies
(\ref{ghat2}).
Moreover, if  $(M^n,g)$ is not conformally diffeomorphic
to the  standard $n-$sphere,  all solutions of (\ref{ghat2}) satisfy
\begin{equation}
\|u\|_{C^{4,\alpha}(M^n,g)}+\|\frac 1u\|_{C^{4,\alpha}(M^n,g)}\le C,
\label{aa7}
\end{equation}
where $C>0$ is some constant depending only on $(M^n,g)$,
 $(f,\Gamma)$ and $\alpha$.
}

\begin{rem} 
$C^0$ and $C^1$ bounds of $u$ and $u^{-1}$ do not depend on the
concavity of $f$. This can be seen from the proof. 
\label{rem1-1}
\end{rem}

 For  $1\le k\le n$, let
$$
\sigma_k(\lambda)=\sum_{1\le i_1<\cdots <i_k\le n}\lambda_{i_1}
\cdots \lambda_{i_k}
$$
be the $k-$th symmetric function and let
$
\Gamma_k$ be the connected component
of $\{\lambda\in \Bbb R^n\ |\ 
\sigma_k(\lambda)>0\}$ containing the positive cone $\Gamma_n$.
 Then, see \cite{CNS},
$(f, \Gamma)=(\sigma_k^{\frac 1k}, \Gamma_k)$ satisfies the hypothesis of 
Theorem~\ref{yamabe}$'$.
\begin{rem} 
For $(f, \Gamma)=(\sigma_1, \Gamma_1)$, 
it is the Yamabe problem in the positive
case  on locally conformally flat
manifolds, and the result is due to
Schoen (\cite{S0}, \cite{S1}).  
For $(f, \Gamma)=(\sigma_2^{\frac 12}, \Gamma_2)$ in
dimension $n=4$, the result was proven without the locally conformally
flatness by Chang, Gursky and Yang (\cite{CGY2}).
For $(f, \Gamma)=(\sigma_n^{ \frac 1n}, \Gamma_n)$, some existence
result was established by Viaclovsky (\cite{V1}) on a class of manifolds.
For  $(f, \Gamma)=(\sigma_k^{\frac 1k}, \Gamma_k)$, the result
was established in our earlier paper \cite{LL}; while the existence
part for $k\ne \frac n2$ was independently established 
by Guan and Wang in \cite{GW2}
using a different method.  Guan, Viaclovsky
and Wang (\cite{GVW}) subsequently proven the algebraic fact that
$\lambda(A_g)\in \Gamma_k$ for $k\ge \frac n2$ implies
the positivity of the Ricci tensor, and therefore $(M,g)$ is
conformally covered by $\Bbb S^n$, and both existence and compactness 
results in this case follow from known results.
For $(f,\Gamma)=(\sigma_k^{\frac 1k}, \Gamma_k)$, $k=3,4$
on $4-$dimensional Riemannian manifolds, as well as for
 $(f,\Gamma)=(\sigma_k^{\frac 1k}, \Gamma_k)$, $k=2,3$, on $3-$dimensional
Riemannian manifolds which are not simply connected , the existence and
compactness results are established by Gursky and Viaclovsky in \cite{GV1}.
\end{rem}
\begin{rem} If we assume in addition that
$f\in C^{k,\alpha}$ for some $k>4$, then, by Schauder theory, 
(\ref{aa7}) can be strengthened as
$$
\|u\|_{C^{k,\alpha}(M^n,g)}+\|u^{-1}\|_{C^{k,\alpha}(M^n,g)}\le C,
$$
where $C>0$ also depends on $k$.
\end{rem}

Since our $C^0$ and $C^1$ estimates for solutions of
(\ref{yamabe1}) (or, equivalently, of   (\ref{ghat2}))
do not make use of the convexity of $V$
(or concavity of $f$), we raise the following

\begin{question}
Under the hypotheses of Theorem \ref{yamabe}$'$,
 but without the
concavity
 assumption on $f$, does there exist a positive Lipschitz function
$u$ on $M^n$ such that $\hat g=u^{\frac 4{n-2}} g$ satisfies
(\ref{ghat2}) in the viscosity sense?
\label{question3}
\end{question}

Equation (\ref{ghat2}) is a fully nonlinear elliptic equation 
of $u$.  Fully nonlinear elliptic equations
 involving $f(\lambda(D^2u))$
have been investigated in the classical and pioneering
paper of Caffarelli, Nirenberg and Spruck \cite{CNS}.  Extensive
studies and outstanding results on such equations
are given by Guan and Spruck \cite{GS}, Trudinger
\cite{T1}, Trudinger and Wang \cite{TW}, and 
many others.
Fully nonlinear equations involving $f(\lambda(\nabla_g^2u+g))$ 
on Riemannian manifolds are studied by Li \cite{Li}, Urbas \cite{U},
and others.   Fully nonlinear equations  involving
the Schouten tensor have been studied by 
Viaclovsky in \cite{V3} and \cite{V1}, and by Chang, Gursky and Yang
in the remarkable papers \cite{CGY2} and \cite{CGY1}.
There have been many papers, preprints,
expository articles, and  works in preparation,  on the
subject and related ones, see, e.g.,
\cite{FG},  \cite{G4},  \cite{V2}  
\cite{GV4}, \cite{GV3}, \cite{GW},
\cite{GW2}, \cite{BGG},
\cite{LL0},  \cite{LL}, \cite{BV}, \cite{GVW}, 
\cite{GV1}, \cite{GV2},  \cite{CY7}, \cite{LL6}, 
\cite{Li5}, \cite{Li8}, \cite{CGY6}, \cite{gw}, \cite{CHY}, \cite{Han},
 \cite{Gon} and \cite{LL5}.
The approach developed in our earlier work \cite{LL} and
continued in  the
present paper makes use of and extends ideas
from previous works on the Yamabe equation by
Gidas, Ni and Nirenberg \cite{GNN},
Caffarelli, Gidas and Spruck \cite{CGS},
Schoen (\cite{S1} and \cite{S2}), Li and
Zhu (\cite{LZhu}), and Li and Zhang (\cite{LZ}).

For $\hat g=u^{\frac 4 {n-2} }g$, we have
(see, e.g., \cite{V3}),
$$
A_{\hat g}=-\frac{2}{n-2}  u^{-1}\nabla^2u+ \frac{2n}{(n-2)^2} u^{-2}
\nabla u \otimes\nabla u -\frac{2}{(n-2)^2} u^{-2}
|\nabla u|^2g+A_{g},
$$
where covariant derivatives on the right side are with respect to $g$.

Let $g_1=u^{\frac{4}{n-2}}g_{flat}$, where $g_{flat}$ denotes
the Euclidean metric on $\Bbb R^n$.
  Then, by the above transformation formula,
$$
A_{g_1}=u^{\frac{4}{n-2}}A^u_{ij}dx^idx^j,
$$
where 
\[
A^u:= -\frac{2}{n-2}u^{  -\frac {n+2}{n-2} }
\nabla^2u+ \frac{2n}{(n-2)^2}u^ { -\frac {2n}{n-2} }
\nabla u\otimes\nabla u-\frac{2}{(n-2)^2} u^ { -\frac {2n}{n-2} }
|\nabla u|^2I,
\]
and  $I$ is the $n\times n$ identity matrix.
In this case, $\lambda(A_{g_1})=\lambda(A^u)$ where $\lambda(A^u)$
denotes the eigenvalues of the $n\times n$ symmetric matrix $A^u$.

Let $\psi$ be a M\"obius transformation
in $\Bbb R^n$, i.e., a transformation
generated by translation, multiplication by nonzero constants, and the
inversion $x \to x/|x|^2$. For any positive $C^2$ function $u$,
let $u_\psi:=|J_\psi|^{ \frac {n-2}{2n} }(u\circ \psi)$ where
$J_\psi$ denotes the Jacobian of $\psi$.  A calculation shows that
$A^{ u_\psi }$ and $A^u\circ \psi$ differ only by
an orthogonal conjugation and therefore
\begin{equation}
\lambda(A^{ u_\psi })=\lambda(A^u)\circ \psi.
\label{conjugate}
\end{equation}

Let ${\cal S}^{n\times n}$ denote the set of
$n\times n$ real symmetric matrices,
 $O(n)$ denote
the set of $n\times n$ real orthogonal matrices,
 $U\subset {\cal S}^{n\times n}$
be  an open set satisfying
\begin{equation}
O^{-1}UO=U,\qquad\forall\ O\in O(n),
\label{U1}
\end{equation}
and let $F\in C^1(U)$ 
 satisfy
\begin{equation}
F(O^{-1}MO)=F(M),\qquad \forall\
M\in U, \ \forall\ O\in O(n),
\label{F1}
\end{equation}

By (\ref{conjugate}) and (\ref{F1}),
\[
F(A^{u_\psi})\equiv F(A^u)\circ \psi.
\]

We proved in \cite{LL} that any conformally invariant operator
$H(\cdot, u, \nabla u, \nabla^2 u)$, in the sense
$$
H(\cdot, u_\psi, \nabla u_\psi, \nabla^2 u_\psi)
\equiv H(\cdot, u, \nabla u, \nabla^2 u)\circ \psi,  
$$
must be of the form $F(A^u)$.

Our next theorem concerns a Harnack type inequality 
for general conformally invariant equations on locally 
conformally flat manifolds.
Let
${\cal S}^{n\times n}_+\subset {\cal S}^{n\times n}$ denote the set
of positive definite matrices.  We will assume that $U$ and $F$ further
satisfy 
\begin{equation}
U\cap \{M+tN\ |\ 0<t<\infty\}\ \mbox{is convex}\qquad \forall\
M\in {\cal S}^{n\times n}, N\in {\cal S}^{n\times n}_+,
\label{U2}
\end{equation}
\begin{equation}
\left(F_{ij}(M)\right)>0,\qquad \forall\ M\in U,
\label{F2}
\end{equation}
where $F_{ij}(M):=\frac{\partial F}{ \partial M_{ij} }(M)$,
and, 
for some $\delta>0$,
\begin{equation}
F(M)\neq 1\qquad\forall\
M\in U\cap \{M\in  {\cal S}^{n\times n}\
|\ \|M\|:=(\sum_{i,j}M_{ij}^2)^{\frac 12} <\delta\}.
\label{F3}
\end{equation}

\begin{thm}
For $n\ge 3$, let $U\subset {\cal S}^{n\times n}$ satisfy
(\ref{U1}) and (\ref{U2}), and let $F\in C^1(U)$ satisfy
(\ref{F1}), (\ref{F2}) and  (\ref{F3}).
For $R>0$, let $u\in C^2(B_{3R})$ be a positive solution of 
\begin{equation}
F(A^u)=1,\quad A^u\in U,\quad \mbox{in}\quad B_{3R},
\label{t2e1}
\end{equation}
where $B_{3R}$ denotes the ball in $\Bbb R^n$ of radius $3R$ and
centered at the origin.
Then
\begin{equation}
(\sup_{B_R}u)(\inf_{B_{2R}}u)\le C(n)\delta^{ \frac{2-n}2 }R^{2-n},
\label{t2e3}
\end{equation}
where $C(n)$ is some constant depending  only on $n$.
\label{t2}
\end{thm}

Let
$$
U_k:=\{ M\in {\cal S}^{ n\times n}\ |\
\lambda(M)\in \Gamma_k\}
$$
and
$$
F_k(M)=\sigma_k(\lambda(M)),\qquad M\in U_k.
$$
For $(F,U)=(F_1, U_1)$, 
(\ref{t2e1}) takes the form
$$
-\Delta u=\frac {n-2}2 u^{ \frac {n+2}{n-2} }, \qquad
\mbox{in}\ B_{3R}.
$$

\begin{rem}  The Harnack type inequality
(\ref{t2e3}) for
$(F,U)=(F_1, U_1)$ was obtained by Schoen in \cite{S2}.
For a class of nonlinearity
including $(F,U)=(F_k^{\frac 1k}, U_k)$, $1\le k\le n$,
 the Harnack type inequality
was established in our earlier work  \cite{LL}.
\end{rem}

\begin{rem} In Theorem \ref{t2},
 there is no concavity assumption
on $F$ and the constant $C(n)$  is given explicitly in the proof.
 The Harnack type inequalities in  \cite{S2}
and \cite{LL} are proved by contradiction arguments which do not yield 
such an explicit constant. 
\end{rem}

Let $g$ be a smooth Riemannian metric on $B_3\subset \Bbb R^n$,
$n\ge 3$, and let
$(f, \Gamma)$ satisfy our usual hypotheses.
Consider
\begin{equation}
f(\lambda(A_{  u^{ \frac 4{n-2} }g }))=1,\quad
\lambda(A_{  u^{ \frac 4{n-2} }g })\in \Gamma, \qquad \mbox{in}\ B_3.
\label{harn1}
\end{equation}

\begin{question}
Are there some positive constants  $C$ and
$\delta$, depending on $(B_3, g)$ and
$(f,\Gamma)$, such that
\[
(\sup_{B_\epsilon} u) 
(\inf_{B_{2\epsilon}} u)\le C \epsilon^{2-n},\qquad
\forall\ 0<\epsilon\le \delta,
\]
holds for any positive solution of (\ref{harn1})?
\label{question4}
\end{question}

\begin{rem}
The answer to the above question is affirmative
for the Yamabe equation
( i.e. $(f,\Gamma)=(\sigma_1, \Gamma_1)$)
in dimension $n=3,4$, see Li and Zhang \cite{LZ2}.
\end{rem}

We have avoided the use of Liouville type theorems
in the proofs of Theorem \ref{yamabe}, Theorem \ref{yamabe}$'$
and Theorem \ref{t2}.  However, in order to  
solve Conjecture \ref{question1} on general Riemannian manifolds,
to answer Question \ref{question4}, 
or to study
 many other issues using fully nonlinear
elliptic equations involving the Schouten tensor, it is important to 
 establish the corresponding Liouville type theorems.

For $n\ge 3$, consider
\begin{equation}
-\Delta u=\frac {n-2}2 u^{ \frac{n+2}{n-2}},
\qquad \mbox{on} \quad \Bbb {R}^n.
\label{eq1new}
\end{equation}

It was proved by Obata (\cite{O}) and Gidas, Ni and Nirenberg (\cite{GNN})
that any  positive $C^2$ solution  of
(\ref{eq1new}) satisfying
$\int_{\Bbb R^n}u^{ \frac {2n}{n-2} }<\infty$
must be of the form
$$
u(x)=(2n)^{ \frac {n-2}4 }
\left(\frac a {1+a^2|x-\bar{x}|^2}\right)^{\frac{n-2}{2}},
$$
where $a>0$ and $\bar x\in \Bbb R^n$.
The hypothesis $\int_{\Bbb R^n}u^{ \frac {2n}{n-2} }<\infty$
was removed by Caffarelli, Gidas and Spruck
(\cite{CGS}); this is important
for applications.
 The method in \cite{GNN} is completely different from that of
\cite{O}.
The method used in
our proof of the Liouville type theorems on general
conformally invariant
fully nonlinear equations 
(Theorem \ref{theorem1}) is in the spirit of
\cite{GNN} rather than that of \cite{O}.
As in \cite{CGS}, the superharmonicity of the solution
has played an important role in our proof of
Theorem \ref{theorem1},  see Lemma \ref{lemma0}.
On the other hand, under some additional hypothesis
on the solution near infinity, the superharmonicity of the solution
is not needed, see theorem 1.4 in \cite{LL}.

Somewhat different proofs of the
 result of
Caffarelli, Gidas and Spruck
 were given
in \cite{CL}, \cite{LZhu} and \cite{LZ}.
In particular, the proofs in \cite{LZhu} and \cite{LZ} 
fully exploit the conformal invariance of the problem
and capture the solutions directly rather
than going through the usual procedure of proving
radial symmetry of solutions and than classifying radial solutions.
A related result of Gidas and Spruck in \cite{GS}
states that there is no positive solution to the
equation $-\Delta u=u^p$ in $\Bbb R^n$
when $1<p<\frac{n+2}{n-2}$.

For $n\ge 3$, $-\infty<p\le\frac{n+2}{n-2}$, we consider the following equation
\begin{equation} 
F(A^u)=u^{p-\frac{n+2}{n-2}},\quad A^u\in U,\quad u>0\quad\mbox{on}~\Bbb{R}^n. 
\label{5} 
\end{equation} 

For $(F,U)=(F_1, U_1)$, equation (\ref{5}) takes the form
$$
-\Delta u=\frac {n-2}2 u^p, \quad u>0, \qquad \mbox{on}\ \Bbb R^n.
$$

\begin{thm} 
For $n\ge 3$, let $U\subset{\cal S}^{n\times n}$ satisfy (\ref{U1}), 
(\ref{U2}), and let $F\in C^1(U)$ satisfy (\ref{F1}), (\ref{F2}). Assume 
that $u\in C^2(\Bbb{R}^n)$ is a superharmonic solution of (\ref{5}) 
for some $-\infty<p\le\frac{n+2}{n-2}$. Then either $u\equiv constant$ 
or $p=\frac{n+2}{n-2}$ and, for some $\bar x\in\Bbb{R}^n$ and some 
positive constants $a$ and $b$ satisfying $2b^2a^{-2}I\in U$ and 
$F(2b^2a^{-2}I)=1$, 
\begin{equation} 
u(x)\equiv (\frac{a}{1+b^2|x-\bar 
x|^2})^{\frac{n-2}{2}},\quad\forall~x\in\Bbb{R}^n.  
\label{6} 
\end{equation} 
\label{theorem1} 
\end{thm} 

\begin{rem} For $(F,U)=(F_k^{\frac 1k}, U_k)$,
$1\le k\le n$,  a solution of
(\ref{5}) is automatically superharmonic.
\end{rem}

\begin{rem} The more difficult case is for  $p=\frac{n+2}{n-2}$.
 When $(F,U)=(F_1, U_1)$, the
result in this case (the rest of this remark also refers to this case),
as mentioned earlier, 
 was established by Caffarelli, Gidas and Spruck (\cite{CGS});
while under some additional hypothesis 
the result was proved by  Obata (\cite{O}) and 
Gidas, Ni and Nirenberg (\cite{GNN}).
For $(F,U)=(F_k^{\frac 1k}, U_k)$, and under some strong hypothesis on $u$ near
infinity, the result was proved by 
Viaclovsky (\cite{V3} and \cite{V2}).  For $
(F,U)=(F_2^{\frac 12}, U_2)$ in dimension $n=4$,
the result was due to Chang,
Gursky and Yang (\cite{CGY2}).  For $(F,U)=(F_k^{\frac 1k}, U_k)$,
 the result was
established in our earlier paper \cite{LL}; while
for $(F,U)=(F_2^{\frac 12}, U_2)$ in dimension $n=5$, as well as
for $(F,U)=(F_2^{\frac 12}, U_2)$
 in dimension $n\ge 6$ under an additional hypothesis
$\int_{\Bbb R^n}u^{ \frac{2n}{n-2} }<\infty$, the result
was independently established by Chang,
Gursky and Yang (\cite{CGY3}).  Under some fairly strong
hypothesis (but weaker
than that used in \cite{V3} and \cite{V2})
on $u$ near infinity, the result was proved in \cite{LL}
without the superharmonicity assumption on $u$.
\end{rem}

\medskip

If we let $(M^n, g)$ denote some smooth compact $n-$dimensional 
Riemannian manifold {\it with boundary}, an analogous problem is to
find conformal metrics with constant scalar
curvature and constant boundary mean curvature.
The problem has been studied by many authors, see, e.g., 
Cherrier (\cite{C}),
Escobar (\cite{E1}, \cite{E3}, \cite{E2} and \cite{E4}),
Han and Li (\cite{HL} and \cite{HL2}), Ambrosetti, Malchiodi
and Li (\cite{AML}), Brendle (\cite{B}), and the references therein.
This boundary Yamabe problem is called of 
positive type if the first eigenvalue of
$$
\left\{
\begin{array}{rll}
-L_g \varphi&=&\lambda \varphi, \qquad 
\mbox{in}\ M^\circ,\\
\frac{\partial \varphi}{\partial \nu}+\frac{n-2}2 h_g\varphi
&=&0,\qquad\qquad \mbox{on}\ \partial M
\end{array}
\right.
$$
is positive, where $h_g$ denotes the mean curvature.

Now we consider an extension of the  
boundary Yamabe problem of positive type to
the fully nonlinear setting:

\begin{question}
Assume that $V$ is an open symmetric convex subset of
$\Bbb {R}^n$, with $\emptyset\neq\partial V\in C^\infty$ satisfying
(\ref{nu1}) and (\ref{nu2}).  Let
$(M^n,g)$ be
 a compact smooth Riemannian manifold with boundary
 satisfying
\[
\lambda(A_g)\in \Gamma(V),\qquad \mbox{on}\ \overline M,
\]
and let $c\in \Bbb R$ be any constant.
Does there exist a smooth positive function
$u\in C^\infty(\overline M)$ such that the conformal metric
 $\hat g=u^{ \frac 4{n-2} }g$ satisfies
(\ref{Q2}) and the boundary mean curvature
$h_{\hat g}$ satisfies
\[
h_{\hat g}=c\qquad\mbox{on}\ \partial M?
\]
\label{question2}
\end{question}

To answer Question \ref{question2}, it is important to
investigate the corresponding Liouville type problem
on half Euclidean space.  Theorem~\ref{liouvilleR+2}
 and Theorem \ref{liouvilleR+3}
below provide such 
Liouville type theorems.

We use $B_R(x)$ to denote the ball in $\Bbb{R}^n$ of
radius $R$ and centered at $x$, and write $B_R=B_R(0)$. Let
$\Bbb{R}^n_{+}=\{(x_1,\cdots,x_n)\in\Bbb{R}^n|~x_n>0\}$ and 
$B^{+}_1=B_1\cap\Bbb{R}^n_{+}$. Consider, for some $c\in \Bbb {R}$,
\begin{equation}
\left\{
\begin{array}{lcl}
F(A^u)=1,\quad A^u\in U,\quad
u>0,\quad&&\mbox{on}~\overline{\Bbb{R}^n_{+}},\\
\frac{\partial u}{\partial
x_n}=cu^{\frac{n}{n-2}},\quad&&\mbox{on}~\partial\Bbb{R}^n_{+}.
\end{array}
\right.
\label{+6a}
\end{equation}

Our first result
is under the assumption that
 the solution has good behavior near infinity.
\begin{thm}
For $n\ge 3$, let $U\subset {\cal S}^{n\times n}$ be an open set
satisfying (\ref{U1}) and (\ref{U2}), and let $F\in C^{1}(U)$
satisfy (\ref{F1}) and (\ref{F2}). For $c\in \Bbb {R}$, we assume that $u\in
C^2(\overline{\Bbb{R}^n_{+}})$ is a solution of (\ref{+6a}) satisfying,
for $u_{0,1}(x):=|x|^{2-n}u(\frac{x}{|x|^2})$,
\begin{equation}
u_{0,1}~~\mbox{can be extended to a positive continuous function
in}~~\overline{B^{+}_1},
\label{+6b}
\end{equation}
\[
\limsup_{x\to 0}(x\cdot\nabla u_{0,1}(x))<\frac{n-2}{2}u_{0,1}(0),
\]
and
\[
\lim_{x\to 0}(|x|^2\nabla u_{0,1}(x))=0.
\]
Then
\begin{equation}
u(x',x_n)\equiv (\frac{a}{1+b|(x',x_n)-(\bar x',\bar
x_n)|^2})^{\frac{n-2}{2}},\quad\mbox{on}~\Bbb{R}^n_{+}, 
\label{11new}
\end{equation}
where $\bar x=(\bar x',\bar x_n)\in\Bbb{R}^n$, $a>0$ and $b+(\min
{\bar x_n,0})^2>0$ are two constants satisfying $2a^{-2}bI\in U$,
$F(2a^{-2}bI)=1$ and $(n-2)a^{-1}b\bar x_n=c$.
\label{liouvilleR+2}
\end{thm}
\begin{rem}
In the above theorem, we do not assume $u$ to be superharmonic.
\end{rem}
\begin{cor} 
For $n\ge 3$, let $U\subset {\cal S}^{n\times n}$ be an open set
satisfying (\ref{U1}) and (\ref{U2}), and let $F\in C^{1}(U)$
satisfy (\ref{F1}) and (\ref{F2}). Assume that $u\in
C^2(\overline{B_1})$ satisfies   
$$
\left\{ 
\begin{array}{lcl} 
F(A^u)=1,~~A^u\in U,~~u>0,\quad&&\mbox{in}~\overline{B_1},\\ 
\frac{\partial u}{\partial 
\nu}+\frac{n-2}{2}u=-cu^{\frac{n}{n-2}},\quad&&\mbox{on}~ \partial B_1, 
\end{array} 
\right. 
$$ 
where $\nu$ denotes the unit outer normal on $\partial B_1$. 
Then $u$ is of the form
\begin{equation}
u(x)\equiv (\frac{a}{1+b|x|^2})^{\frac{n-2}{2}}\quad\mbox{in}~B_1,
\label{+6}
\end{equation}
where $a,b,c$ satisfy
\[ 
a>0,\quad \frac{2b}{a^2}I\in U,\quad 
F(\frac{2b}{a^2}I)=1,\quad\frac{n-2}{2}(1-b)=-ca. 
\]  
\label{corA1} 
\end{cor}

Our next Liouville type theorem does not require any hypothesis
on the solution near infinity.
\begin{thm} 
For $n\ge 3$, let $U\subset {\cal S}^{n\times n}$ be an open set
satisfying (\ref{U1}) and (\ref{U2}), and let $F\in C^{1}(U)$
satisfy (\ref{F1}) and (\ref{F2}). Assume that 
\begin{equation} 
0\notin\overline {F^{-1}(1)}. 
\label{more1} 
\end{equation} 
For $c\in \Bbb {R}$, we assume that $u\in 
C^2(\overline{\Bbb{R}^n_{+}})$ is a solution of (\ref{+6a}) satisfying, 
\[
\Delta u\le 0\quad\mbox{in}~\Bbb{R}^n_{+}. 
\] 
Then $u$ is of form (\ref{11new}) with $\bar x,a$ and $b$ given below 
(\ref{11new}).  
\label{liouvilleR+3} 
\end{thm} 
\begin{rem} 
For $c\le 0$, the assumption (\ref{more1}) is not needed. This can be 
seen in the proof.  
\end{rem} 
\begin{rem} 
$(F,U)=(F_k^{\frac 1k}, U_k)$,
$1\le k\le n$,  satisfy the hypotheses of the theorem.    
\end{rem}
\begin{rem}
For $(F,U)=(F_1, U_1)$, the result was proved by Li and Zhu 
\cite{LZhu}; while under an additional hypothesis
$u(x)=O(|x|^{2-n})$ for large $|x|$, the solutions were
classified by Escobar \cite{E1}.
\end{rem} 

Our proofs of 
  Theorem~\ref{liouvilleR+2} and Theorem \ref{liouvilleR+3} make use
of the following result concerning radially symmetric solutions.
\begin{thm}
For $n\ge 3$, let $U\subset {\cal S}^{n\times n}$ be an open set
satisfying (\ref{U1}), and let $F\in C^{1}(U)$
satisfy (\ref{F1}) and (\ref{F2}). Assume that $u\in C^2(B_1)$ is
radially symmetric and satisfies
\[
F(A^u)=1,\quad A^u\in U,\quad u>0,\quad\mbox{in}~B_1.
\]
Then $u$ is of the form (\ref{+6}) with $a>0$, $b\ge -1$, $\frac{2b}{a^2}I\in
U$ and $F(\frac{2b}{a^2}I)=1$. 
\label{liouvilleR+1}
\end{thm}

In the following we state some
of the results in a forthcoming paper
\cite{LL5}.  
First, an  existence and compactness
result on subcritical equations:
\begin{thm}
Let $(M,g)$ be a smooth, compact, connected Riemannian manifold of
dimension $n\ge 3$, and let $1<1+\epsilon\le p\le
\frac{n+2}{n-2}-\epsilon<\frac{n+2}{n-2}$. Then there exists a
positive solution $u\in C^\infty(M)$ to 
\begin{equation}
\sigma_k^{\frac 1k}(A_{u^{\frac
{4}{n-2}}g})=u^{p-\frac{n+2}{n-2}},\quad\mbox{on}~M.
\label{subyamabe}
\end{equation}
Moreover all positive solutions of (\ref{subyamabe}) satisfy, 
 for all $m\ge 2$,
\[
\|u\|_{C^{m}(M,g)}+\|\frac 1 u\|_{C^{m}(M,g)}\le C,
\] 
where $C>0$ depends only on $(M^n,g)$, $\epsilon$ and $m$. 
\label{pcompact}
\end{thm}

\begin{rem}  For $k=1$, this is well known.
\end{rem}

Next, a Harnack type inequality on 
half Euclidean balls:
\begin{thm}
For $n\ge 3$ and $R>0$, let $u\in C^2(\overline {B^+_{3R}})$ be a
solution of the equation  
\[
\left\{
\begin{array}{lcl}
\sigma_k^{\frac 1k}(A^u)=1,&&\quad\mbox{in}~B^+_{3R}:=B_{3R} \cap\Bbb
{R}^n_+.\\ 
\frac{\partial u}{\partial x_n}=cu^{\frac
{n}{n-2}},&&\quad\mbox{on}~\partial B_{3R}^+\cap\partial\Bbb
{R}^n_+\quad\mbox{for some constant}~c.\\ 
u>0,~A^u\in\Gamma_k,&&\quad\mbox{on}~\overline {B^+_{3R}}.
\end{array}
\right.
\]
Then there exists some constant $C>0$ depending only on $n$ and $c$
such that
\[
(\sup_{B_R^+}u)(\inf_{\partial B_{2R}^+}u)\le CR^{2-n}.
\]
\label{harnack+}
\end{thm}

\begin{rem}
For $k=1$, this, as well as a stronger form, is
established by Li and Zhang in \cite{LZ} (see theorem 1.7 and remark 1.11
there).
\end{rem}

\begin{rem}
Theorem~\ref{pcompact} and Theorem~\ref{harnack+} hold for more
general $(f,\Gamma)$ (see \cite{LL5}).
\end{rem}

As mentioned earlier, Theorem \ref{yamabe}$'$
in the case $(f, \Gamma)=(\sigma_1, \Gamma_1)$ is the Yamabe problem
in the positive case
on locally conformally flat manifolds,
 and the result is due to Schoen
(\cite{S0} and \cite{S1}).
  The proof in \cite{S1}
has three main ingredients: The first is the
 existence of the developing map due to Schoen and Yau \cite{SY1},
 the second is 
the use of  the method of moving planes, and the third is
the Liouville type theorem of 
Caffarelli, Gidas and Spruck \cite{CGS}.
A major difficulty in extending the result for $(f, \Gamma)
=(\sigma_1, \Gamma_1)$
to fully nonlinear $(f,\Gamma)$ was the lack of 
corresponding Liouville type theorem.
An important step was taken by Zhang and the second author
in \cite{LZ} which 
gives a proof of Schoen$'$s Harnack type inequality
for the Yamabe equation without using the Liouville type theorem
in \cite{CGS}.
Adapting this idea, we established in  \cite{LL} (theorem 1.27 there)
the Harnack type inequality (\ref{t2e3}) for
a class of nonlinearity
including $(F,U)=(F_k^{\frac 1k}, U_k)$, $1\le k\le n$, under
the circumstance that the corresponding Liouville type theorem was
not available. This also made us  recognize 
the possibility of proving Theorem \ref{yamabe}$'$ 
without the corresponding Liouville type theorem.
Indeed we have developed in \cite{LL}
an approach, based on the method of moving spheres
(i.e. the method of moving planes together
with the conformal invariance of the problem),  
to prove the existence and compactness
results for the fully  
nonlinear version of the Yamabe problem on 
locally conformally flat manifolds under the circumstance that
the corresponding Liouville type theorem was not available.
Another major difficulty in proving 
Theorem \ref{yamabe}$'$ is the lack of $C^0$ and $C^1$
estimates of solutions.  
We have developed a new approach
in \cite{LL}, again  based on the method of moving 
spheres,  to obtain such estimates. 
We have also introduced in \cite{LL} a homotopy
which connects the general
fully nonlinear version of the Yamabe
problem to the Yamabe problem   and used the degree for
second order fully nonlinear elliptic operators in \cite{Li1}
and the result in \cite{S1} for the Yamabe problem
to prove the existence of solutions to the fully nonlinear ones.

In \cite{GLW} Guan, Lin and Wang have also presented a proof of
Theorem \ref{t2} under
an additional concavity hypothesis
on $F$ and of  Theorem \ref{yamabe}$'$.  
We clarify these overlaps in this paragraph:
First,  these results
follow immediately from our earlier work \cite{LL}
and Lemma \ref{l2} ---- a quantitative version of
a calculus lemma used repeatedly in \cite{LL}.
 Second  we completed the proof of these results earlier.
Indeed, the only change
one  needs to make is to  move the four lines below (4.3) 
on page 1446 of \cite{LL} to be right after line 5 of the same page.
After making this change, the
 gradient estimate stated on line 7 of the same page
follows from Lemma \ref{l2}, and
Theorem \ref{t2} under
an additional concavity hypothesis
on $F$ and Theorem \ref{yamabe}$'$, as well as our
new $C^0$ and $C^1$ estimates, follow from  the proofs of 
theorem 1.25 and theorem 1.27 in \cite{LL}.
We did not see the elementary proof of Lemma \ref{l2}
at the time of submitting \cite{LL} to the journal, but proved it
soon afterwards.  Theorem \ref{t2}
and  Theorem \ref{yamabe}$'$, with an emphasis on our
new $C^0$ and $C^1$ estimates based on the method of moving planes,
were presented by the second author in his $45-$minute
invited talk at ICM 2002 in August 2002 in Beijing.
Told us by C.S. Lin that he
 started to work with G. Wang
in October-November 2002 which
led to \cite{GLW} where 
a proof of Theorem \ref{t2} under
an additional concavity hypothesis
on $F$ and   Theorem \ref{yamabe}$'$ is included.
The proof,  following \cite{LL}
(in particular following the above mentioned steps developed
there), provides the only ingredient beyond \cite{LL}
  which,
 as explained above, 
 amounts to the calculus lemma (Lemma \ref{l2}).
We present the proof
of  Theorem \ref{yamabe}$'$ and  Theorem \ref{t2} 
  in Section 2 and Section 3 respectively.  The proof
of  Theorem \ref{yamabe}$'$,
appeared in slightly shorter form in \cite{LL6}
and  in preprint form  \cite{LL7}, contains one slight
simplification to the arguments in \cite{LL} which 
avoids the use of local $C^2$ estimates (only global $C^2$ estimates
are needed); while the
 proof of
 Theorem \ref{t2}, also appeared in slightly shorter form in \cite{LL6}
and in   \cite{LL7},
 contains one more ingredient to remove the concavity
assumption on $F$ which also yields an explicit constant $C(n)$
in (\ref{t2e3}).

Due to Theorem \ref{yamabe}$'$ (or Theorem \ref{yamabe}),
 Conjecture
 \ref{question1}$'$ (or Conjecture \ref{question1})
mainly concerns the problem on Riemannian manifolds which are not
locally conformally flat.  In general, 
Equation (\ref{ghat2}) does not have a variational formulation.
A plausible  approach is to establish a priori 
estimates (\ref{aa7}) for all solutions of (\ref{ghat2}), and 
to use the homotopy in \cite{LL} to connect the problem to the Yamabe
problem.   For the Yamabe problem (i.e. (\ref{ghat2}) for
$(f,\Gamma)=(\sigma_1, \Gamma_1)$ ), such estimate
was given by Li and Zhu \cite{LZhu} in dimension $n=3$;
the estimate in dimension $n=4$  follows from  a combination 
of the results of Li and Zhang \cite{LZ2} and
Druet \cite{D1}; Li and Zhang have 
extended the estimate to dimension $n\le 7$, 
as well as to dimension $n\ge 8$ but
under an additional hypothesis that the
Weyl tensor of $g$ is nowhere vanishing, 
 see \cite{LZ3}.
The Liouville type theorem of Caffarelli, Gidas and Spruck
has played an important role in the proof of this result.
It is clear that Theorem \ref{theorem1} will also play
an important role in proving Conjecture \ref{question1}$'$.

The main difficulty in proving Theorem \ref{theorem1}
is to remove the possible isolated singularity of $u$ at infinity.
By the conformal invariance of the problem, 
we may assume that the isolated singularity is
at $0$ instead of at infinity.
The following analytical issue is relevant:  Let
$u\in C^\infty(B_1\setminus \{0\})$ and $v\in C^\infty(B_1)$
be positive solutions of
\[
F(A^u)=1, \quad A^u\in U,\qquad
\mbox{in} \ B_1\setminus \{0\},
\]
and 
\[
F(A^v)=1, \quad A^v\in U,\qquad
\mbox{in} \ B_1,
\]
satisfying
\[
u>v\qquad\mbox{in}\ B_1\setminus \{0\}.
\]
Is it true that 
\[
\liminf_{|x|\to 0}(u(x)-v(x))>0?
\]

If the answer to the above question  were ``yes'', then
the proof of theorem 1.4 in \cite{LL} 
would yield a proof of   
Theorem \ref{theorem1} for $p=\frac{n+2}{n-2}$.
So far, the answer to the question 
 is not known even for 
$(F,U)=(F_k^{\frac 1k}, U_k)$, $2\le k\le n$.
The answer to the question is ``yes'' for
$(F,U)=(F_1, U_1)$ due to  some 
elementary properties of superharmonic functions in a punctured ball.
As far as we know, 
the isolated singularity issue 
encountered in the application of the method of moving plane
has always been handled by providing
an affirmative answer to a local question like the above.  
Our proof of Theorem \ref{theorem1} avoids 
this local question by exploiting global information
of $u$, through a delicate use of Lemma \ref{lemma0}.
The proof of Theorem \ref{theorem1} also fully 
exploits the conformal invariance of the problem
and captures the solutions directly rather
then going through the usual procedure of proving
radial symmetry of solutions and than classifying radial solutions. 
Two proofs of  Theorem \ref{theorem1}
appeared in preprint forms in \cite{LL3} and \cite{LL4}.
We present in Section 4 the proof in \cite{LL4}.
Theorem \ref{liouvilleR+2} and Theorem \ref{liouvilleR+3},
some Liouville type theorems on half Euclidean spaces,
are extensions of theorem 1.4
in \cite{LL} and Theorem \ref{theorem1} respectively.
The proofs are given in Section 5.

\bigskip

\noindent{\bf Acknowledgment.}\
Part of this paper was completed while the second author 
was a visiting member at the Institute for Advanced Study
in Fall 2003.
He thanks J. Bourgain and IAS for providing him the 
excellent environment, as well as for providing him the financial
support through NSF-DMS-0111298.
Part of the work of the second author is also supported by
NSF-DMS-0100819.

\section{\bf Proof of Theorem \ref{yamabe} and Theorem \ref{yamabe}$'$.}
In Appendix B, we deduce the equivalence of Theorem~\ref{yamabe} and
Theorem~\ref{yamabe}$^{\prime}$., therefore we only need to prove one
of the two theorems.\newline
\noindent {\bf Proof of Theorem~\ref{yamabe}$^{\prime}$.}\ Without loss of
generality, we further assume $f$ is homogeneous  of 
degree $1$. Indeed, in Appendix B, we construct a new
function $\tilde f$ which is homogeneous of degree $1$, satisfies
the same assumptions as $f$ does, and $\tilde f^{-1}(1)=
f^{-1}(1)$.

We first establish (\ref{aa7}). Let $(\widetilde M,
\widetilde g)$ be the universal cover of $(M^n,g)$,
 with $i:\widetilde M\to M^n$ being a covering map and $ \widetilde g=i^*g$. 
It is well-known that there exists a conformal immersion
$$
\Phi: (\widetilde M,  \widetilde g) \to (\Bbb S^n, g_0),
$$
where $ g_0$ denotes the standard metric on $\Bbb S^n$. 
By $\lambda (A_g)\in\Gamma$ and the assumption $\Gamma\subset
 \Gamma_1$, we have $R_g> 0$. Hence by a deep theorem of Schoen and
Yau in \cite{SY1}, $\Phi$
 is injective. Let  
$$
\Omega=\Phi(\widetilde M).
$$
\begin{clm}
\[
\frac 1C\le u\le C,\quad |\nabla_g u|\le C\qquad\mbox{on}\ M^n,
\]
where $u\in C^2(M^n)$ is an arbitrary positive 
solution of (\ref{ghat2}) with $\hat g=u^{  \frac 4{n-2} }g$ and $C>0$
is some constant depending only on $(M^n,g)$ and $(f,\Gamma)$.
\label{estimateclaim}
\end{clm}

For convenience, we introduce
$$
U=\{A\in {\cal S}^{n\times n}\ |\ 
\lambda(A)\in \Gamma\},
$$
and 
$$
F(A)=f(\lambda(A)),\qquad A\in U.
$$

We distinguish two cases.

\noindent {\bf Case 1.}\
$$ 
\Omega= \Bbb S^n;
$$
\noindent {\bf Case 2.}\ 
$$
\Omega\neq \Bbb S^n.
$$

 In Case 1,
$(\Phi^{-1})^*\widetilde g=\eta^{\frac 4{n-2}}g_0$ on
$\Bbb S^n$, where $\eta$ is a positive smooth function on
$\Bbb S^n$. 
Let $\tilde u= u\circ i$.  Since  
$F\left(A_{ \tilde u^{\frac 4{n-2}}
\widetilde g }
\right) =1$ on $\widetilde M$, we have 
$$
F\left(
A_{ [(\tilde u\circ \Phi^{-1})\eta]^{\frac 4{n-2}}g_0 }\right)=1,
\qquad \mbox{on}\  \Bbb S^n.
$$ 

By corollary~1.6 in \cite {LL}, 
$ (\tilde u\circ \Phi^{-1})\eta=a |J_{ \varphi}|^{ \frac {n-2}{2n} }$
for some positive constant $a$ and some conformal diffeomorphism
$\varphi: \Bbb S^n\to \Bbb S^n$.  Since
$\varphi^* g_0=|J_{ \varphi}|^{ \frac 2n}g_0$,
we have, by the above equation, that
$$
f(a^{-\frac 4{n-2}}(n-1)e)=
f(a^{-\frac 4{n-2}}\lambda(A_{g_0}))=1,
$$
where $e=(1,\cdots,1)$. By (\ref{aa6}) and the concavity of $f$, we
know $\nabla f(\lambda)\cdot\lambda>0$ for any $\lambda\in\Gamma$. Thus
$f|_{\partial\Gamma}=0$ and (\ref{aa6}) implies $a$ is a constant
uniquely determined by $(f,\Gamma)$.

Fix a compact subset $E$ of $\widetilde M$ such that
$i(E)=M^n$.
Since  $(M^n,g)$ is not conformally
diffeomorphic to $(\Bbb S^n, g_0)$,
so $\pi_1(M^n)$ is nontrivial.  Let $\tilde x^{(1)}\in E$ and
$\tilde x^{(2)}\in \widetilde M$ 
be two distinct points satisfying
$\tilde u(\tilde x^{(1)})=\tilde u(\tilde x^{(2)})=\max\limits_{M^n}u$.
Then 
$$
dist_{g_0}\left( \Phi(\tilde x^{(1)}),  \Phi(\tilde x^{(2})\right) 
\ge \frac 1C.
$$
Consequently,
$$
\min\{ |J_\varphi(\Phi(\tilde x^{(1)})|,  |J_\varphi(\Phi(\tilde x^{(2})| \}
\le C,
$$
from which we deduce that
$$
\min\{ \tilde u(\tilde x^{(1)})\eta (\Phi(\tilde x^{(1)})),
 \tilde u(\tilde x^{(2)})\eta (\Phi(\tilde x^{(2)})) \}\le C.
$$
It follows that
$$
\max\limits_{M^n}u=\tilde u(\tilde x^{(1)})=\tilde u(\tilde x^{(2)})
\le C.
$$
Moreover, we also know from the above and the formula of
$\tilde u$ that
$$
|J_\varphi|\le C\qquad\mbox{on}\ \Bbb S^n,
$$
from which we deduce that 
$$
\||J_\varphi|\|_{ C^m(\Bbb S^n, g_0) }
+\|\frac 1{ |J_\varphi| }\|_{ C^m(\Bbb S^n, g_0) } \le C(m)
$$
and therefore
$$
\|u\|_{ C^m(M^n,g) }+\|u^{-1}\|_{ C^m(M^n,g) } \le C
$$
for some $C$ depending only on $(M,g)$, $(f,\Gamma)$ and
$m$. Estimates (\ref{aa7}) is established in this case.

In Case 2,  by the result in \cite{SY1}, $\Omega=\Phi(\tilde{M})$ is an
open and dense subset of $\Bbb S^n$, $(\Phi^{-1})^*\widetilde
g=\eta^{\frac 4{n-2}}g_0$ on $\Omega$, where $\eta$ is a positive
smooth function in $\Omega$ satisfying
$\lim\limits_{z\to \partial\Omega}\eta(z)=\infty$.
Let $u(x)=\max\limits_{M^n} u$ for some $x\in M^n$, and let
$i(\tilde x)=x$ for some $\tilde x\in E$.
By composing with a rotation of $\Bbb S^n$, we
may assume without loss of generality that  
$\Phi(\tilde x)=S$, the south pole of $\Bbb S^n$.
Let $P: \Bbb S^n \to \Bbb R^n$ be the stereographic projection, and
let $v$ be  the positive function on the open subset
 $P(\Omega)$ of $\Bbb R^n$
determined by
$(P^{-1})^*(\eta^{ \frac 4{n-2} }g_0)=v^{ \frac 4{n-2} }g_{flat}$,
where $g_{flat}$ denotes the
Euclidean metric on $\Bbb R^n$. 
Then for some $\epsilon>0$, depending only on $(M^n,g)$, 
we have
$$
B_{9\epsilon}:=\{x\in \Bbb R^n\
|\ |x|<9\epsilon\}\subset P(\Omega),
$$
and
$$
dist_{ flat} \big(P(\Phi(E)),
\partial  P(\Omega)\big)>9\epsilon.
$$

On  $P(\Omega)$, 
$$
F(A^{\hat u})=1,\qquad \lambda(A^{\hat u})\in \Gamma,
$$
where $\hat u=(\tilde u\circ \Phi^{-1}\circ P^{-1})v$.

By the property of $\eta$, we know that
\begin{equation}
\lim\limits_{P(\Omega)\ni y\to\bar y\in\partial P(\Omega)}\hat u(y)=\infty,
\label{unbounded}
\end{equation}
and, if the north pole of $S^n$ does not belong to $\Omega$,
\begin{equation}
\lim\limits_{y\in P(\Omega),|y|\to\infty}(|y|^{n-2}\hat u (y))=\infty.
\label{49-1}
\end{equation}

For every $x\in \Bbb R^n$ satisfying 
$dist_{flat}(x, P(\Phi(E)))<2\epsilon$, we can 
perform a moving sphere argument as in the
corresponding part in  \cite{LL}  (for $w_j$ there)
to show that, $\forall\ 0<\lambda<4\epsilon,
|y-x|\ge\lambda,~y\in P(\Omega)$, 
\begin{equation}
\hat u_{x,\lambda}(y):=\frac { \lambda^{n-2} }
{ |y-x|^{n-2} } \hat u( \frac { \lambda^2 (y-x) }{ |y-x|^2 } )\le
\hat u(y).
\label{move}
\end{equation}

When proving the above, there is 
some minor difference between the north pole of $S^n$, $N\in\Omega$
and $N\notin\Omega$. If $N\notin\Omega$, then by (\ref{49-1}), there
is no worry about ``touching at infinity'' in the moving sphere
procedure.  If $N\in\Omega$, then $\infty$ is a regular point of $\hat u$ 
(i.e., $|z|^{2-n}\hat u(\frac z{|z|^2})$ can be extended as a $C^2$
positive function near $z=0$) and therefore  by the strong maximum
principle argument as in \cite{LL}, if ``touching at infinity'' occurs, 
$(\hat u)_{x,\lambda}$ would coincide with $\hat u$ in the unbounded connected
component of $P(\Omega)$  for some 
$0<\lambda<4\epsilon$, which violates (\ref{unbounded}) since
$(\hat u)_{x,\lambda}$ is apparently bounded near any
point of $\partial P(\Omega)$.

By Lemma~\ref{l2} in Appendix A, we deduce from (\ref{move}) that 
\[
|\nabla(log \hat u)(y)|\le C(\epsilon)\quad
\forall\ dist_{flat}(y, P(\Phi(E)))< \epsilon.
\]
It follows, for some $C$ depending only on $(M^n,g)$, that
$$
|\nabla_g \log u|\le C\qquad \mbox{on}\ M^n.
$$
Hence Claim~(\ref{estimateclaim}) follows directly from the bounds below
\begin {equation}
\min\limits_{M^n} u\le C,\quad \max\limits_{M^n} u\ge \frac 1C\quad\mbox{for
some universal constant}~C.
\label {VC}
\end{equation}
To establish (\ref{VC}), let $u(\bar x)=\min\limits_{M^n} u$. At $\bar
x$, by $\nabla u (\bar x)=0$, $(\nabla ^2 u (\bar x))\ge 0$, and
(\ref{aa5}), we have  
\[
1=f(\lambda
(A_{\hat g}))\le f(u^{\frac{-4}{n-2}}\lambda (A_g)),
\] 
which implies, by $f|_{\partial \Gamma}=0$ and $f\in C^0(\bar\Gamma)$,
that $u^{\frac{-4}{n-2}}(\bar x)\ge C$, i.e., $u(\bar x)\le
C$. Similarly, by properties of $f$ (in particular (\ref{aa6})), we
can establish $\max\limits_{M^n} u\ge \frac 1C$.
The $C^2$ estimate of $u$ has been established in
\cite{LL} (see also \cite {V1} for the estimates for
$(f,\Gamma)=(\sigma_k^{\frac 1k},\Gamma_k)$). The $C^2$ estimate of
$u^{-1}$ follows in view of Claim~\ref{estimateclaim}.

Thus when $(M^n,g)$ is not conformally diffeomorphic to a standard sphere, we
have proved that any positive solution of (\ref{ghat2}) satisfies, for
some constant $C$ depending only on $(M^n,g)$ and $(f,\Gamma)$, 
$$
\|u\|_{ C^2(M^n,g)}+\|u^{-1}\|_{ C^2(M^n,g)}\le C.
$$

Since $f$ is concave in $\Gamma$, $C^{2,\alpha}$ and higher order
derivative estimates follow from a theorem of
Evans (\cite {E}) and Krylov (\cite{K}), and the Schauder
estimate. 

To establish the existence part of Theorem~\ref{yamabe}$^{\prime}$, 
we only need to treat the case that $(M^n,g)$ is not conformally
diffeomorphic to a standard sphere since it is obvious otherwise.
We use the following homotopy introduced in \cite{LL}. For
$0\le t\le 1$,  let
\[
f_t(\lambda)=f\left(t\lambda+(1-t)\sigma_1(\lambda)e\right),
\]
be defined on
\[
\Gamma_t:=\{ \lambda\in \Bbb R^n\ |\
t\lambda+(1-t)\sigma_1(\lambda)e\in \Gamma\},
\]
where $e=(1,1,\cdots,1)$.

Consider, for $0\le t\le 1$,
\begin{equation}
f_t(\lambda(A_{\hat g}))=1,\quad\lambda(A_{\hat g})\in \Gamma_t,
\qquad \mbox{on}\ M^n.
\label{equationt}
\end{equation}
Here and below $\hat g=u^{ \frac 4{n-2}}g$.

By the a priori estimates we have just established, there exists  some
constant $C>0$  independent of $t\in [0,1]$ such that for all solutions
$u$ of (\ref{equationt}),
\begin{equation}
\|u\|_{C^{4,\alpha}}(M^n,g)+\|u^{-1}\|_{C^{4,\alpha}}(M^n,g)\le C.
\label{uu1}
\end{equation}

By (\ref{uu1}) and the assumption $f|_{\partial\Gamma}=0$, $\exists~\delta>0$
independent of $t\in [0,1]$  such that all solutions $u$ of
(\ref{equationt}) satisfy 
$$
dist(\lambda(A_{\hat g}), \partial \Gamma_t)\ge 2 \delta.
$$

Define, for $0\le t\le 1$,
\begin{eqnarray*}
O_t^*=\{u\in C^{4,\alpha}(M^n)\ &|&\
\lambda(A_{\hat g})\in \Gamma_t,\ dist(\lambda(A_{\hat g}),
 \partial \Gamma_t)> \delta,\
\\
&&\qquad u>0,\
 \|u\|_{C^{4,\alpha}(M^n,g)}+\|u^{-1}\|_{C^{4,\alpha}(M^n,g)}<2C\},
\end{eqnarray*}
where $C$ is the constant in (\ref{uu1}).
By \cite{Li1},
$$
d_t:=\deg \left( F_t-1, O_t^*,0\right),\qquad
0\le t\le 1,
$$
is well defined, where $F_t[u]:=f_t(\lambda(A_{\hat g}))-1$, and 
$$
d_t\equiv d_0,\qquad 0\le t\le 1.
$$
In particular,
$$
d_1=d_0.
$$
The equation (\ref{equationt}) for $t=0$ is the Yamabe equation. By
the result of Schoen in \cite{S1} for the Yamabe problem, $d_0=-1$.
Thus $d_1\neq 0$ and equation (\ref{ghat2}) has a solution. Theorem
\ref{yamabe}$^{\prime}$ is established.

\section{Proof of Theorem~\ref{t2}}
\noindent{\bf Proof of Theorem~\ref{t2}.}\  
Part of the proof of this theorem is taken from 
\cite{LL}, which we include here for reader$'$s convenience.
We only need to prove the theorem for
$R=\delta=1$.  Indeed, let
$$
\tilde F(M):=F(\delta M),
\quad \tilde U:=\delta^{-1}U, \quad
\mbox{and}\ \ \tilde u(x):=\delta^{ \frac { n-2} 4 }
R^{  \frac { n-2} 2} u(Rx).
$$
Then 
$$
\tilde F(A^{\tilde u})=1,\quad  A^{\tilde u}\in 
\tilde U,\quad \mbox{in}\quad B_{3},
$$
and $(\tilde F, \tilde U)$ satisfies the hypothesis
of Theorem \ref{t2} with $R=\delta=1$. Thus, once we
have established the theorem in the case $R=\delta=1$,
we have
$$
(\sup_{B_R}u)(\inf_{B_{2R}}u)=
\delta^{ \frac{2-n}2 }R^{2-n}
(\sup_{B_1}\tilde u)(\inf_{B_2}\tilde u)
\le C\delta^{ \frac{2-n}2 }R^{2-n}.
$$

Thus we assume in the following $R=\delta=1$.
Let $u(\bar{x})=\max\limits_{\bar{B}_1}u$. 
As in the proof of theorem 1.27 in \cite{LL},we can find $\tilde{x}\in
B_{\frac 12}(\bar{x})$ such that
\[
u(\tilde{x})\ge 2^{\frac{2-n}{2}}\sup_{B_{\sigma}(\tilde{x})}u
\]
and
\begin{equation}
\gamma:=u(\tilde{x})^{\frac{2}{n-2}}\sigma\ge\frac1 2
u(\bar{x})^{\frac{2}{n-2}},
\label{another2}
\end{equation}
where $\sigma=\frac12(1-|\tilde{x}-\bar{x}|)\le\frac12$.

If 
$$
\gamma\le 2^{n+8}n^4,
$$
then
$$
(\sup_{B_1}u)(\inf_{B_2}u)
\le u(\bar x)^2\le (2\gamma)^{ \frac {n-2}2 }
\le C(n),
$$
and we are done.  So we always assume that
$$
\gamma>  2^{n+8}n^4.
$$

Let $\Gamma:=u(\tilde{x})^{\frac{2}{n-2}}\ge 2\gamma$, and consider 
\[
w(y):=\frac{1}{u(\tilde{x})}u\Big(\tilde{x}+\frac{y}{u(\tilde{x})
^{\frac{2}{n-2}}}\Big),\quad |y|<\Gamma.
\]
Clearly
\begin{equation}
\min_{\partial B_\Gamma}
w\ge \frac{1}{u(\tilde{x})}\inf_{B_2}u,
\label{another3}
\end{equation}
\begin{equation}
1=w(0)\ge 2^{\frac{2-n}{2}} \sup_{ B_\gamma}w.
\label{another4}
\end{equation}
By the conformal invariance of the equation satisfied by $u$,
$$
F(A^w)=1, \qquad w>0,\qquad \mbox{on}\ B_\Gamma.
$$

Fix
\[
r= 2^{n+6} n^4<\frac 1{4} \gamma.
\]
$\forall |x|<r$, 
 consider 
$$
w_{x,\lambda}(y):= (\frac \lambda{ |y-x| })^{n-2}
w(x+ \frac {\lambda^2(y-x) }{ |y-x|^2 }).
$$
By the conformal invariance of the equation, we have
$$
F(A^{ w_{x,\lambda}})=1, \qquad w_{x,\lambda}>0,\qquad \mbox{on}\ 
B_\Gamma \setminus B_\lambda(x),\quad\forall~0<\lambda<\frac{3\gamma}{4}.
$$ 
As in \cite{LL},  there exists $0<\lambda_x<r$ such that
we have
$$
w_{x,\lambda}(y)\le w(y),\qquad
\forall \ 0<\lambda<\lambda_x, \ y\in B_\Gamma\setminus B_\lambda(x),
$$
and
$$
w_{x,\lambda}(y)< w(y),\qquad
\forall \ 0<\lambda<\lambda_x, \ y\in \partial  B_\Gamma.
$$
By the moving sphere argument as in \cite{LL}, 
we only need to consider the following two cases:

\noindent {\bf Case 1.}\  For some $|x|<r$ and some
 $\lambda\in (0,r)$,   $w_{x,\lambda}$ touches
$w$ on $\partial B_\Gamma$.

\noindent {\bf Case 2.}\  For all $|x|<r$ and all 
$\lambda\in (0,r)$, we have 
$$
w_{x,\lambda}(y)\le w(y),\quad \forall\  
|y-x|\ge\lambda, \ y\in B_\Gamma.
$$

In  Case 1, let $\lambda\in (0,r)$ be the smallest
number for which  $w_{x,\lambda}$ touches $w$
on $\partial B_\Gamma$.
By (\ref{another3}), we have, for some $|y_0|=\Gamma$,
\[
\frac{1}{u(\tilde{x})}\inf_{B_2}u\le\min_{\partial B_\Gamma}w=
w_{x,\lambda}(y_0). 
\]
Recall (\ref{another4}), 
\[
w_{x,\lambda}(y_0)\le \Big(\frac{\lambda}{|y_0-x|}\Big)^{n-2}
\sup_{B_\gamma}w\le 
2^{\frac{n-2}{2}}\Big(\frac{\lambda}{|y_0-x|}\Big)^{n-2}\le
2^{\frac{n-2}{2}}\Big(\frac{r}{\Gamma-r}\Big)^{n-2}.
\]
Therefore
\[
\sigma^{\frac{n-2}{2}}u(\tilde{x})\inf_{B_2}u\le
2^{\frac{n-2}{2}}\sigma^{\frac{n-2}{2}}u(\tilde{x})^2
\Big(\frac{r}{\Gamma-r}\Big)^{n-2}.
\]
Since $4r<\gamma\le\frac{\Gamma}{2}$ and $\sigma\le\frac12$,
\begin{equation}
\sigma^{\frac{n-2}{2}}u(\tilde{x})\inf_{B_2}u\le
2^{\frac{n-2}{2}}\sigma^{\frac{n-2}{2}}u(\tilde{x})^2\frac{r^{n-2}}
{(\frac12\Gamma)^{n-2}}=
2^{\frac32(n-2)}\sigma^{\frac{n-2}{2}}r^{n-2}
\le 2^{n-2}r^{n-2}.
\label{another5}
\end{equation}
We deduce from  (\ref{another2}) and (\ref{another5})  that
\[
(\sup_{B_1}u)(\inf_{B_2}u)\le 4^{n-2}r^{n-2}
\le C(n).
\]

In Case 2, we have, by Lemma~\ref{l2}
and (\ref{another4}), that  
$$
|\nabla w(y)|\le 2(n-2)r^{-1}w(y)\le
(n-2)2^{ \frac {n}2 }r^{-1},\quad\forall |y|\le r.
$$
Let $\epsilon$ be the number such that
$$
\xi(y):=\frac{1-\epsilon}{r}(r-|y|^2),
\qquad |y|<\sqrt{r}
$$
satisfies
\[
w\ge \xi,\quad\mbox{on}~B_{\sqrt{r}},
\]
and, for some $|\bar y|<\sqrt{r}$,
$$
w(\bar y)=\xi(\bar y).
$$
Since $1=w(0)\ge \xi(0)=1-\epsilon$ and 
$w(\bar y)>0$, we have $0\le \epsilon<1$.

By the estimates of $|\nabla w|$ and the mean value theorem,
$$
|w(y)-1|=|w(y)-w(0)|\le (n-2)2^{ \frac {n}2} r^{-\frac 12},
\qquad \forall\ |y|\le \sqrt{r}.
$$
So 
$$
1-  (n-2)2^{ \frac {n}2} r^{-\frac 12}\le w(\bar y)=\xi(\bar y)\le 
1-\epsilon,
$$
and therefore
$$
0\le \epsilon\le  (n-2)2^{ \frac {n}2} r^{-\frac 12}.
$$
Clearly,
$$
\nabla w(\bar y)=\nabla \xi(\bar y),
\ \ |\nabla \xi(\bar y)|\le \frac 2{\sqrt{r}},
\ \ D^2w(\bar y)\ge D^2\xi(\bar y)= -2(1-\epsilon)r^{-1}I.
$$
It follows that
\[
A^w(\bar y)\le
A^{\xi}(\bar y)\le
\frac{(10n+4)}{ (n-2)^2 } 2^{ \frac {2n}{n-2} } r^{-1}I.
\]
Since $F(A^w(\bar y))=1$, we have, by (\ref{F3}) (recall
that $\delta=1$),
$ \frac{(10n+4)}{ (n-2)^2 } 2^{ \frac {2n}{n-2} } r^{-1}
\ge 1$, violating the choice of $r$.
Thus we have shown that Case 2 can never occur.
Theorem~\ref{t2} is established.

\vskip 5pt
\hfill $\Box$
\vskip 5pt

\section{\bf Proof of Theorem~\ref{theorem1}}
\begin{lem}
For $n\ge 2$, $B_1\subset\Bbb {R}^n$, let $u\in L^1_{loc}(B_1\setminus
\{0\})$ be the solution of  
\[
\Delta u\le 0\qquad\mbox{in}~B_1\setminus \{0\}
\]
in the distribution sense. Assume $\exists$ $a\in R$ and $p\neq q\in\Bbb {R}^n$
such that
\[
u(x)\ge \max \{a+p\cdot x-\delta (x),a+q\cdot x-\delta
(x)\}\quad\forall x\in~B_1\setminus \{0\},
\]
where $\delta(x)\ge 0$ satisfies $\lim\limits_{x\to 0}\frac {\delta
(x)}{|x|}=0$. Then
$$
\lim\limits_{r\to 0}\inf\limits_{B_r}u>a.
$$
\label{lemma0}
\end{lem}
\noindent{\bf Proof.}\ Let 
\[
v(x):=a+p\cdot x-\delta (x),\quad w(x):=a+q\cdot x-\delta
(x),\quad\forall x\in B_1.
\]
By subtracting
$a+p\cdot x$ from
$u,v,w$ respectively, we can assume $a=0$ and $p=0$. After a
rotation and a dilation of the coordinates, we can also assume $\nabla
w(0)=e_1$.

Let $u_{\epsilon}:=\frac{1}{\epsilon}u(\epsilon\cdot)$,
$v_{\epsilon}:=\frac{1}{\epsilon}v(\epsilon\cdot)$, and
$w_{\epsilon}:=\frac{1}{\epsilon}w(\epsilon\cdot)$. We have 
\[
v_{\epsilon}(x)=o(1),\quad w_{\epsilon}(x)=x_1+o(1),
\]  
where $o(1)\to 0$ uniformly on $\bar {B}_1$ as $\epsilon\to 0$.
$\forall \bar \delta>0$, by $u_{\epsilon}\ge v_{\epsilon}$, $\exists
\epsilon_0>0$ such that
\[
u_{\epsilon}(x)\ge -\bar \delta\quad\mbox{in}~
B_1,\quad\forall~\epsilon\le\epsilon_0.  
\]
By $u_{\epsilon}\ge
w_{\epsilon}$, we have $u_{\epsilon}\ge c_0>0$ on $\Omega:=B_{\frac
14}(\frac 12 e_1)$ for some universal constant $c_0$ independent of
$\bar \delta$ and $\epsilon$.

Let $\xi^{\bar\delta}$ be the solution of
\[
\left\{
\begin{array}{lcl}
\Delta\xi^{\bar\delta}=0\qquad&&\mbox{in}~B_1\setminus\bar\Omega \\
\xi^{\bar\delta}=\frac{c_0}{2}~~\mbox{on}~\partial\Omega,\quad
\xi^{\bar\delta}=-2\bar\delta~~&&\mbox{on}~\partial B_1.
\end{array}
\right.
\]
Since $\xi^{\bar\delta}\to \xi^0$ in $C^{\infty}(\bar {B}_1)$, we
have, for small $\bar\delta$,
\begin{equation}
\xi^{\bar\delta}(0)>\frac 12 \xi^0(0)>0,
\label{positive}
\end{equation}
where $\xi^0$ is the solution of 
\[
\left\{
\begin{array}{lcl}
\Delta\xi^0=0\qquad\mbox{in}~B_1\setminus\bar\Omega&&\\
\xi^0=\frac{c_0}{2}~~\mbox{on}~\partial\Omega,\quad
\xi^0=0~~\mbox{on}~\partial B_1.&& 
\end{array}
\right.
\]
In the following, we fix some $\bar\delta>0$ such that (\ref{positive})
holds.
 
Let $G$ be the solution of
$$
\left\{
\begin{array}{lcl}
-\Delta G=\delta_0\qquad\mbox{in}~B_1&&\\
G=0~~\mbox{on}~\partial B_1,\quad
G(x)\to\infty~~\mbox{as}~x\to 0,&& 
\end{array}
\right.
$$
where $\delta_0$ is the Dirac mass at $0$.

Let $A>1$ be chosen later. $\forall 0<\delta<\frac {1}{10}$, consider
$\eta_{\epsilon}:=u_{\epsilon}+\frac{A}{G(\delta)}G-\xi^{\bar\delta}$ on
$B_1\setminus \{B_{\delta}\cup\Omega\}$. We have
$$
\Delta\eta_{\epsilon}\le 0\quad\mbox{in}~B_1\setminus\overline
{\{B_{\delta}\cup\Omega\}}.  
$$
Near $\partial B_{\delta}$,
\[
\eta_{\epsilon}\ge -\bar\delta+\frac A2-\frac 12c_0>0,\quad\mbox{for
large}~A,
\]
and near $\partial B_1$,
$\eta_{\epsilon}\ge-\bar\delta+\frac 32\bar\delta>0$. Hence 
\begin{equation}
\eta_{\epsilon}>0\quad\mbox{in}~B_1\setminus\overline
{\{B_{\delta}\cup\Omega\}}. 
\label{eta}
\end{equation}
For any fixed $x\in B_1\setminus\{0\}$, $\forall~0<\delta<|x|$,
$\forall \epsilon>0$ small, sending $\delta\to 0$ in (\ref{eta}), it
leads to $u_{\epsilon}(x)\ge\xi^{\bar\delta}(x)$. Therefore,
$\forall\epsilon\le \epsilon_0$,
\[
\lim\limits_{r\to 0}\inf\limits_{B_r}u=\lim\limits_{r\to
0}\inf\limits_{B_r}u_{\epsilon}\ge \xi^{\bar\delta}(0)>\frac 12\xi_0(0)>0.
\]  
\vskip 5pt
\hfill $\Box$
\vskip 5pt
Lemma~\ref{lemma0} is sufficient for our use. Such result holds for
more general linear elliptic operators of second
order. For example, we have 
\begin{lem}
For $n\ge 2$ and $B_1\subset\Bbb {R}^n$, let $u\in C^2(B_1\setminus
\{0\})$ satisfy  
\[
Lu:=a^{ij}u_{ij}+b^iu_i+cu\le f\qquad\mbox{in}~B_1\setminus \{0\},
\]
where $(a^{ij})>0$ and $a^{ij}\in
C^{\alpha}(B_1)$ for some $0<\alpha<1$, $f,b_i,c\in
L^{\infty}(B_1)$. Assume $\exists$ $a\in R$ and $p\neq q\in\Bbb {R}^n$
such that
\[
u(x)\ge \max \{a+p\cdot x-\delta (x),a+q\cdot x-\delta
(x)\}\quad\forall x\in~B_1\setminus \{0\},
\]
where $\delta>0$ satisfies $\lim\limits_{x\to 0}\frac {\delta (x)}{|x|}=0$.
Then
\[
\liminf\limits_{x\to 0}u(x)>a.
\]
\label{ge1l2}
\end{lem}
\noindent{\bf Proof of Lemma~\ref{ge1l2}.}\ Let 
\[
v(x):=a+p\cdot x-\delta (x),\quad w(x):=a+q\cdot x-\delta
(x),\quad\forall x\in B_1.
\]
By subtracting
$a+p\cdot x$ from $u,v,w$ respectively, and replacing $f(x)$
by $f(x)-b^i(x)v_i(0)-c(x)\nabla v(0)\cdot x-c(x)v(0)$, we can assume
\[
a=0,\quad p=0,\quad Lu\le f\quad\mbox{in}~B_1\setminus\{0\}. 
\] 
Let $Q\in GL(n)$ satisfy $Q(a^{ij}(0))Q^t=I_{n\times n}$. Replacing
$u,v,w$ by 
\[
u(Q^{-1}\cdot),\quad v(Q^{-1}\cdot),\quad w(Q^{-1}\cdot),
\] 
and $a^{ij}(x),b^i(x),c(x),f(x)$ by 
\[
Q(a^{ij}(Q^{-1}x))Q^t,\quad Q^t(b^i(Q^{-1}x)),\quad c(Q^{-1}x),\quad
\bar {f} (Q^{-1}x) 
\] 
respectively, we can assume $(a^{ij})(0)=I_{n\times n}$.

Let $u_{\epsilon}:=\frac{1}{\epsilon}u(\epsilon\cdot)$,
$v_{\epsilon}:=\frac{1}{\epsilon} v(\epsilon\cdot)$, and
$w_{\epsilon}:=\frac{1}{\epsilon} w(\epsilon\cdot)$. We have 
\[
v_{\epsilon}(x)=o(1),\quad w_{\epsilon}(x)=\nabla w(0)\cdot
x+o(1)\quad\mbox{on}~ \bar B_1. 
\]  
We may also assume that $|\nabla w(0)|=1$ by a dilation. Hence, by
$u_{\epsilon}\ge v_{\epsilon}$ and $u_{\epsilon}\ge w_{\epsilon}$,
$\forall \bar\delta>0$, $\exists \epsilon_0>0$ such that $\forall
\epsilon\le \epsilon_0$,  
\[
u_{\epsilon}(x)\ge -\bar\delta~~\mbox{on}~ B_1,\quad
u_{\epsilon}\ge c_0~~\mbox{on}~ \Omega:=B_{\frac 14}(\frac 12 \nabla w(0)),
\]
where $c_0>0$ is some universal constant independent of $\bar\delta$
and $\epsilon$. Moreover $u_{\epsilon}$ satisfies the equation
\[
L^{\epsilon}u_{\epsilon}(x):= a^{ij}(\epsilon
x)(u_{\epsilon})_{ij}(x)+\epsilon b^i(\epsilon
x)(u_{\epsilon})_i(x)+\epsilon^2 c(\epsilon x) u_{\epsilon}(x)\le
\epsilon f(\epsilon x)\quad\mbox{in}~B_1. 
\]

Let $\xi_{\bar\delta}$ be the solution of 
$$
\left\{
\begin{array}{lcl}
L^{\epsilon}\xi^{\bar\delta}(x)=\epsilon f(\epsilon
x)\qquad\mbox{in}~B_1\setminus\bar\Omega&&\\   
\xi^{\bar\delta}=\frac{c_0}{2}~~\mbox{on}~\partial\Omega,\quad
\xi^{\bar\delta}=-2\bar\delta~~\mbox{on}~\partial B_1.&& 
\end{array}
\right.
$$
We have $\xi^{\bar\delta}\to\xi^0$ in $C^1(B_1\setminus\bar\Omega)$,
where $\xi^0$ is the solution of 
$$
\left\{
\begin{array}{lcl}
\Delta \xi^0=0\qquad\mbox{in}~B_1\setminus\bar\Omega&&\\ 
\xi^0=\frac{c_0}{2}~~\mbox{on}~\partial\Omega,\quad
\xi^0=0~~\mbox{on}~\partial B_1.&& 
\end{array}
\right.
$$
Hence we can initially pick some $\bar\delta>0$ such that
$\xi^{\bar\delta}(0)>\xi^0(0)>0$. 

Let $G$ be the solution of
$$
\left\{
\begin{array}{lcl}
-L^{\epsilon} G=\delta_0\qquad\mbox{in}~B_1&&\\
G=0~~\mbox{on}~\partial B_1,\quad
G(x)\to\infty~~\mbox{as}~x\to 0.&& 
\end{array}
\right.
$$
We know $G$ is asymptotically radial as $\epsilon\to 0$.

Let $A>1$ be chosen later. $\forall 0<\delta<\frac {1}{10}$, consider
$\eta_{\epsilon}:=u_{\epsilon}+\frac{A}{\min\limits_{\partial
B_{\delta}}G}G-\xi^{\bar\delta}$ on $B_1\setminus
\{B_{\delta}\cup\Omega\}$. We have 
$$
L^{\epsilon}\eta_{\epsilon}\le 0\quad\mbox{in}~B_1\setminus\overline
{\{B_{\delta}\cup\Omega\}}.  
$$
On $\partial B_{\delta}$,
\[
\eta_{\epsilon}\ge -\bar\delta+A-\frac{c_0}{2}>0,
\]
and on $\partial B_1$, $-\bar\delta+2\bar\delta=\bar\delta>0$. Hence 
\begin{equation}
\eta_{\epsilon}>0\quad\mbox{in}~B_1\setminus\overline
{\{B_{\delta}\cup\Omega\}}. 
\label{eta1}
\end{equation}
For any fixed $x\in B_1\setminus\{0\}$, $\forall~0<\delta<|x|$,
$\forall \epsilon\le\epsilon_0$, sending $\delta\to 0$ in
(\ref{eta1}), then $u_{\epsilon}(x)\ge\xi^{\bar\delta}(x)$. Therefore
$\liminf\limits_{x\to 0}u(x)=\liminf\limits_{x\to 0}u_{\epsilon}
(x)\ge \xi^{\bar\delta}(0)>\frac 12\xi^0(0)>0$.

\vskip 5pt
\hfill $\Box$
\vskip 5pt

\noindent 
{\bf Proof of Theorem~\ref{theorem1} for $p=\frac{n+2}{n-2}$.}\ Since 
$u$ is a positive superharmonic function, we have, by the maximum 
principle, that 
\[ 
u(x)\ge\frac{\min\limits_{\partial B_1}u}{|x|^{n-2}},\quad\forall~|x|\ge 1. 
\] 
In particular  
\begin{equation} 
\liminf\limits_{|x|\to\infty}(|x|^{n-2}u(x))>0. 
\label{7} 
\end{equation} 
\begin{lem} 
For any $x\in\Bbb{R}^n$, there exists $\lambda_0(x)>0$ such that  
\[ 
u_{x,\lambda}(y):=(\frac{\lambda}{|y-x|})^{n-2} 
u(x+\frac{\lambda^2(y-x)}{|y-x|^2})\le 
u(y),\quad\forall~|y-x|\ge\lambda,~0<\lambda<\lambda_0(x).  
\] 
\label{lemma1} 
\end{lem} 
{\bf Proof.}\ This follows from the proof of 
lemma~2.1 in \cite{LZ}.

\vskip 5pt 
\hfill $\Box$ 
\vskip 5pt 
 
For any $x\in\Bbb{R}^n$, set  
\[ 
\bar\lambda(x):=\sup\{\mu~|~u_{x,\lambda}(y)\le 
u(y),~\forall~|y-x|\ge\lambda,~0<\lambda<\mu\}.  
\] 
Let  
\begin{equation} 
\alpha:=\liminf\limits_{|x|\to\infty}(|x|^{n-2}u(x)). 
\label{8} 
\end{equation} 
Because of (\ref{7}),  
\begin{equation} 
0<\alpha\le\infty. 
\label{8a} 
\end{equation} 
If $\alpha=\infty$, then the moving sphere procedure will never stop 
and therefore $\bar\lambda(x)=\infty$ for any $x\in\Bbb{R}^n$. This 
follows from arguments in \cite{LZ} and \cite{LL} (see also
\cite{LL2}). By the definition of $\bar\lambda(x)$ and the fact
$\bar\lambda(x)=\infty$, we have,   
\[ 
u_{x,\lambda}(y)\le u(y),\quad\forall~|y-x|\ge\lambda>0. 
\] 
By a calculus lemma (see e.g., lemma~11.2 in \cite{LZ}), $u\equiv 
constant$, and Theorem~\ref{theorem1} for $p=\frac{n+2}{n-2}$ is 
proved in this case (i.e. $\alpha=\infty$). So, from now on, we assume 
 
\begin{equation} 
0<\alpha<\infty. 
\label{9} 
\end{equation} 
By the definition of $\bar \lambda(x)$, 
\[ 
u_{x,\lambda}(y)\le 
u(y),\quad\forall~|y-x|\ge\lambda,~0<\lambda<\bar\lambda (x). 
\] 
Multiplying the above by $|y|^{n-2}$ and sending $|y|\to\infty$, we 
have, 
\[ 
\alpha\ge\lambda^{n-2}u(x),\qquad\forall\ 0<\lambda<\bar\lambda(x). 
\] 
Sending $\lambda\to\bar\lambda(x)$, we have (using (\ref{9})), 
\begin{equation} 
\infty>\alpha\ge\bar\lambda(x)^{n-2}u(x),\quad\forall~x\in\Bbb{R}^n. 
\label{10} 
\end{equation} 
Since the moving sphere procedure stops at $\bar\lambda(x)$, we must 
have, by using the arguments in \cite{LZ} and \cite{LL} (see also
\cite{LL2}),   
\begin{equation} 
\liminf\limits_{|y|\to\infty}(u(y)-u_{x,\bar\lambda(x)}(y))|y|^{n-2}=0, 
\label{11} 
\end{equation} 
i.e., 
\begin{equation} 
\alpha=\bar\lambda(x)^{n-2}u(x),\quad\forall~x\in\Bbb{R}^n. 
\label{12} 
\end{equation} 
Let us switch to some more convenient notations. For a Mobius 
transformation $\phi$, we use notation  
\[ 
u_{\phi}:=|J_{\phi}|^{\frac{n-2}{2n}}(u\circ\phi), 
\] 
where $J_{\phi}$ denotes the Jacobian of $\phi$.

For $x\in\Bbb{R}^n$, let 
\[ 
\phi^{(x)}(y):=x+\frac{\bar\lambda(x)^2(y-x)}{|y-x|^2}, 
\]  
we know that 
$u_{\phi^{(x)}}=u_{x,\bar\lambda(x)}$.

Let $\psi(y):=\frac{y}{|y|^2}$, and let 
\[ 
w^{(x)}:=(u_{\phi^{(x)}})_{\psi}=u_{\phi^{(x)}\circ\psi}. 
\]  
For $x\in\Bbb{R}^n$, the only possible singularity for $w^{(x)}$ (on 
$\Bbb{R}^n\cup\{\infty\}$) is $\frac{x}{|x|^2}$. In particular, $y=0$ 
is a regular point of $w^{(x)}$. A direct calculation yields 
\[ 
w^{(x)}(0)=\bar\lambda(x)^{n-2}u(x), 
\] 
and therefore, by (\ref{12}), 
\[ 
w^{(x)}(0)=\alpha,\quad\forall~x\in\Bbb{R}^n. 
\]
Clearly, $u_{\psi}\in C^2(\Bbb{R}^n\setminus\{0\})$ and $\Delta
u_{\psi}\le 0$ in $\Bbb{R}^n\setminus\{0\}$, $\liminf\limits_{y\to
0}u_{\psi}(y)=\alpha$, and, for some $\delta(x)>0$,  
\[ 
w^{(x)}\in C^2(B_{\delta (x)}),\quad\forall~x\in\Bbb{R}^n, 
\] 
\[ 
u_{\psi}\ge 
w^{(x)}\quad\mbox{in}~B_{\delta(x)}\setminus\{0\},\quad\forall~x\in\Bbb{R}^n.  
\]  
\begin{lem} 
$\nabla w^{(x)}(0)=\nabla w^{(0)}(0)$, i.e., $\nabla w^{(x)}(0)$ is 
independent of $x\in\Bbb{R}^n$. 
\label{lemma2} 
\end{lem} 
{\bf Proof of Lemma~\ref{lemma2}.}\ This follows from Lemma 
\ref{lemma0}. Indeed, 
for any $x,~\tilde x\in\Bbb{R}^n$, let 
\[ 
v:=w^{(x)},\quad w:=w^{(\tilde x)},\quad u:=u_{\psi}. 
\] 
We know that $w(0)=v(0)$, $u_{\psi}\ge w$ and $u_{\psi}\ge v$ 
near the origin, and we also know that 
 $\liminf\limits_{y\to 0}u_{\psi}(y)=w(0)$, so, 
by Lemma \ref{lemma0},  we must have $\nabla 
v(0)=\nabla w(0)$, i.e., $\nabla w^{(x)}(0)=\nabla w^{(\tilde 
x)}(0)$. Lemma~\ref{lemma2} is established. 
 
\vskip 5pt 
\hfill $\Box$ 
\vskip 5pt

For $x\in\Bbb{R}^n$, 
\begin{eqnarray*} 
w^{(x)}(y)&=&\frac{1}{|y|^{n-2}} 
\Big\{(\frac{\bar\lambda(x)}{|\frac{y}{|y|^2}-x|})^{n-2}u(x+ 
\frac{\bar\lambda(x)^2(\frac{y}{|y|^2}-x)}{|\frac{y}{|y|^2}-x|^2})\Big\}\\ 
&=&(\frac{\bar\lambda(x)}{|\frac{y}{|y|}-|y|x|})^{n-2}u(x+ 
\frac{\bar\lambda(x)^2(y-|y|^2 x)}{|\frac{y}{|y|}-|y|x|^2})\\ 
&=&(\frac{\bar\lambda(x)^2}{1-2x\cdot y+|y|^2 |x|^2})^{\frac{n-2}{2}}u(x+ 
\frac{\bar\lambda(x)^2(y-|y|^2x)}{1-2x\cdot y+|y|^2|x|^2}). 
\end{eqnarray*} 
So, for $|y|$ small, 
\[ 
w^{(x)}(y)=\bar\lambda(x)^{n-2}(1+(n-2)x\cdot y)u(x+\bar\lambda(x)^2 y)+O(|y|^2), 
\] 
and, using (\ref{12}), 
\[ 
\nabla 
w^{(x)}(0)=(n-2)\bar\lambda(x)^{n-2}u(x)x+\bar\lambda(x)^n\nabla 
u(x)=(n-2)\alpha x+\alpha^{\frac{n}{n-2}}u(x)^{\frac{n}{2-n}}\nabla u(x). 
\] 
By Lemma~\ref{lemma2}, $\vec V:=\nabla w^{(x)}(0)$ is a constant 
vector in $\Bbb{R}^n$, so we have, 
\[ 
\nabla_x(\frac{n-2}{2}\alpha^{\frac{n}{n-2}}u(x)^{-\frac{2}{n-2}} 
-\frac{(n-2)\alpha}{2}|x|^2+\vec V\cdot x)\equiv 0. 
\] 
Consequently, for some $\bar x\in\Bbb{R}^n$ and $d\in R$,  
\[ 
u(x)^{-\frac{2}{n-2}}\equiv\alpha^{-\frac{2}{n-2}}|x-\bar 
x|^2+d\alpha^{-\frac{2}{n-2}}.  
\] 
Since $u>0$, we must have $d>0$. Thus 
\[ 
u(x)\equiv (\frac{\alpha^{\frac{2}{n-2}}}{d+|x-\bar x|^2})^{\frac{n-2}{2}}. 
\] 
Let $a=\alpha^{\frac{2}{n-2}}d^{-1}$ and $b=d^{-\frac 12}$. Then $u$ 
is of the form (\ref{6}). Clearly $A^u(0)=2b^2a^{-2}I$, so 
$2b^2a^{-2}I\in U$ and $F(2b^2a^{-2}I)=1$. Theorem~\ref{theorem1} in 
the case $p=\frac{n+2}{n-2}$ is established.

\vskip 5pt 
\hfill $\Box$ 
\vskip 5pt 
 
{\bf Proof of Theorem~\ref{theorem1} for 
$-\infty<p<\frac{n+2}{n-2}$.}\   In this case, the equation satisfied by 
$u$ is no longer conformally invariant, but it transforms to our 
advantage when making reflections with respect to spheres, i.e., the 
inequalities have the right direction so that the strong maximum 
principle and the Hopf lemma can still be applied. 

First, we still have (\ref{7}) since this only requires the 
superharmonicity and the positivity of $u$. Lemma~\ref{lemma1} still 
holds since it only uses (\ref{7}) and the $C^1$ regularity of $u$ in 
$\Bbb{R}^n$. For $x\in\Bbb{R}^n$, we still define $\bar\lambda(x)$ in 
the same way. We also define $\alpha$ as in (\ref{8}) and we still 
have (\ref{8a}). 

For $x\in\Bbb{R}^n$, $\lambda>0$, the equation of $u_{x,\lambda}$ now 
takes the form 
\begin{equation} 
F(A^{u_{x,\lambda}}(y))=(\frac{\lambda}{|y-x|})^{(n-2)(\frac{n+2}{n-2}-p)} 
u_{x,\lambda}(y)^{p-\frac{n+2}{n-2}},\quad A^{u_{x,\lambda}}(y)\in 
U,\quad\forall~y\neq x. 
\label{14} 
\end{equation} 
\begin{lem} 
If $\alpha=\infty$, then $\bar\lambda(x)=\infty$ for all $x\in\Bbb{R}^n$. 
\label{lemma3} 
\end{lem} 
{\bf Proof.}\ Suppose the contrary, 
$\bar\lambda(\bar x)<\infty$ for some $\bar x\in\Bbb{R}^n$. Without 
loss of generality, we may assume $\bar x=0$, and we use notations  
\[ 
\bar\lambda:=\bar\lambda(0),\quad u_{\lambda}:=u_{0,\lambda},\quad 
B_{\lambda}:= B_{\lambda}(0). 
\]  
By the definition of $\bar\lambda$,  
\[
u_{\bar\lambda}\le u\quad\mbox{on}~~\Bbb{R}^n\setminus B_{\bar\lambda}. 
\] 
By (\ref{14}), 
\begin{equation} 
F(A^{u_{\bar\lambda}})\le u_{\bar\lambda}^{p-\frac{n+2}{n-2}},\quad 
A^{u_{\bar\lambda}}\in U,\quad\mbox{on}~\Bbb{R}^n\setminus B_{\bar\lambda}. 
\label{16} 
\end{equation} 
Recall that $u$ satisfies  
\begin{equation} 
F(A^u)=u^{p-\frac{n+2}{n-2}},\quad A^u\in 
U,\quad\mbox{on}~~\Bbb{R}^n\setminus B_{\bar\lambda}. 
\label{17} 
\end{equation} 
By (\ref{16}) and (\ref{17}), 
\begin{equation} 
F(A^{u_{\bar\lambda}})-F(A^u) 
-(u_{\bar\lambda}^{p-\frac{n+2}{n-2}}-u^{p-\frac{n+2}{n-2}})\le 
0,\quad A^{u_{\bar\lambda}}\in U,~A^u\in 
U,\quad\mbox{on}~~\Bbb{R}^n\setminus B_{\bar\lambda}. 
\label{18} 
\end{equation} 
Since $\alpha=\infty$, we have  
\begin{equation} 
\liminf\limits_{|y|\to\infty}|y|^{n-2}(u-u_{\bar\lambda})(y)>0. 
\label{19} 
\end{equation} 
The inequality in (\ref{18}) goes the right direction. Thus, with 
(\ref{19}), the arguments for $p=\frac{n+2}{n-2}$ work essentially in 
the same way here and we obtain a contradiction by continuing the 
moving sphere procedure a little bit further. This deserves some 
explanations. Because of (\ref{19}), and using arguments in \cite{LL}
(see also \cite{LL2}), we only need to show that 
\begin{equation} 
u_{\bar\lambda}(y)< u(y),\quad\forall~|y|>\bar\lambda, 
\label{19a} 
\end{equation} 
and 
\begin{equation} 
\frac{d}{dr}(u-u_{\bar\lambda})|_{\partial B_{\bar\lambda}}>0, 
\label{19b} 
\end{equation} 
where $\frac{d}{dr}$ denotes the differentiation in the outer normal 
direction with respect to $\partial B_{\bar\lambda}$.

If $u_{\bar\lambda}(\bar y)= u(\bar y)$ for some $|\bar 
y|>\bar\lambda$, then, using (\ref{18}) as in the proof of lemma~2.1 
in \cite{LL}, we know that $u_{\bar\lambda}-u$ satisfies that 
\[ 
L(u_{\bar\lambda}-u)\le 0, 
\] 
where $L=-a_{ij}(x)\partial_{ij}+b_i(x)\partial_i+c(x)$ with 
$(a_{ij})>0$ continuous and $b_i$, $c$ continuous. 

Since $u_{\bar\lambda}-u\le 0$ near $\bar y$, we have, by the strong 
maximum principle, $u_{\bar\lambda}\equiv u$ near $\bar y$. For the 
same reason, $u_{\bar\lambda}(y)\equiv u(y)$ for any 
$|y|\ge\bar\lambda$, violating (\ref{19}). (\ref{19a}) has been 
checked. Estimate (\ref{19b}) can be established in a similar way by 
using the Hopf lemma (see the proof of lemma~2.1 in \cite{LL}). Thus 
Lemma~\ref{lemma3} is established. 
 
\vskip 5pt 
\hfill $\Box$ 
\vskip 5pt

By Lemma~\ref{lemma3} and the usual arguments, we know that if 
$\alpha=\infty$, $u$ must be a constant, and Theorem~\ref{theorem1} 
for $-\infty<p<\frac{n+2}{n-2}$ is also proved in this case. 

From now on, we always assume (\ref{9}). As before, we obtain 
(\ref{10}). Since the inequality in (\ref{16}) goes the right 
direction, the arguments for $p=\frac{n+2}{n-2}$ (see also the 
arguments in the proof of Lemma~\ref{lemma3}) essentially apply and we 
still have (\ref{11}) and (\ref{12}). Applying the rest of the arguments for 
$p=\frac{n+2}{n-2}$, we have $u$ is of the form 
(\ref{6}) with some positive constants $a$ and $b$. However, we know 
that, for $u$ of the form (\ref{6}), $A^u\equiv 2b^2a^{-2}I$ and 
$F(A^u)\equiv constant$. This violates (\ref{5}) since 
$u^{p-\frac{n+2}{n-2}}$ is not a constant when
$p<\frac{n+2}{n-2}$. Theorem~\ref{theorem1} for
$-\infty<p<\frac{n+2}{n-2}$ is established.   
 
\vskip 5pt 
\hfill $\Box$ 
\vskip 5pt

\section{\bf Proof of Theorem~\ref{liouvilleR+2}-\ref{liouvilleR+1}}
\subsection{\bf Proof of Theorem~\ref{liouvilleR+2}.}
To prove Theorem~\ref{liouvilleR+2}, let us first establish
Theorem~\ref{liouvilleR+1}.

\noindent  
{\bf Proof of Theorem~\ref{liouvilleR+1}.}\ Let $u$ be the same
as in Theorem~\ref{liouvilleR+1}, we let 
\[
v(r)=u(r,0\cdots,0),\quad 0\le r<1.
\]
Clearly, $v'(0)=0$. For $x=(r,0\cdots,0)$, $0<r<1$, we have
\[
\nabla u(x)=(v'(r),0\cdots,0),\quad\nabla^2
u(x)=diag(v''(r),\frac{v'(r)}{r},\cdots, \frac{v'(r)}{r}),
\]
and 
\[
A^u(x)=diag(\lambda_1^v(r),\lambda_2^v(r),\cdots,\lambda_n^v(r)),
\]
where 
$$
\left\{
\begin{array}{lcl}
&&\lambda_1^v(r)=-\frac{2}{n-2}v^{-\frac{n+2}{n-2}}v''
+\frac{2(n-1)}{(n-2)^2}v^{-\frac{2n}{n-2}}(v')^2\\
&&\lambda_2^v(r)=\cdots=\lambda_n^v(r)=
-\frac{2}{n-2}v^{-\frac{n+2}{n-2}}\frac{v'}{r}
-\frac{2}{(n-2)^2}v^{-\frac{2n}{n-2}}(v')^2.
\end{array}
\right.
$$
Here and in the following, we use $diag(\lambda_1,\cdots,\lambda_n)$
to denote the diagonal matrix 
\[
\left(
\begin{array}{cccc}
\lambda_1&&&\\
&\lambda_2&&\\
&&\ddots&\\
&&&\lambda_n
\end{array}
\right)
\]
Let $w(x)=(\frac{a}{1+b|x|^2})^{\frac{n-2}{2}}$ with
$a=v(0)^{\frac{2}{n-2}}$ and
$b=\frac{1}{2-n}a^{\frac{2-n}{2}}v''(0)$. With these choices of
$a$ and $b$, we have 
\[
w(0)=v(0),\quad w'(0)=v'(0)=0,\quad w''(0)=v''(0).
\]
A calculation yields 
\[
A^w(x)\equiv\frac{2b}{a^2}I=A^u(0),
\]
and therefore $w$ satisfies 
\[
F(A^w)=1,\quad A^w\in U,\quad w>0,\quad\mbox{in}~\{x\in\Bbb{R}^n~|~b|x|^2>-1\}.
\]
Introduce
$f(\lambda_1,\cdots,\lambda_n)=F(diag(\lambda_1,\cdots,\lambda_n))$. Clearly,
\[
\lambda_j(0):=\lim\limits_{r\to
0}\lambda_j(r)=-\frac{2}{n-2}v(0)^{-\frac{n+2}{n-2}}v''(0),\quad 1\le
j\le n,
\]
and therefore, by the symmetry of $f$ in $\lambda_1\cdots,\lambda_n$,
we have
\[
f_{\lambda_j}(\lambda_1(0),\cdots,\lambda_n(0))
=f_{\lambda_1}(\lambda_1(0),\cdots,\lambda_n(0)),\quad 2\le j\le n.
\]
Since $diag(\lambda_1(0),\cdots,\lambda_n(0))\in U$, we have, by
(\ref{F2}), $f_{\lambda_1}(\lambda_1(0),\cdots,\lambda_n(0))>0$. 
\begin{lem}
Let $\alpha$ and $\beta$ be positive constants, and let $k\ge 1$ be an
integer satisfying $k+\gamma>\alpha$ for some $0<\gamma\le 1$. Assume
that $\xi\in C^{k-1,\gamma}([0,\beta])$ satisfies  
\begin{equation}
|\xi(r)|\le \frac{\alpha}{r}\int_0^r |\xi (s)|~ds,\quad\forall~0<r<\beta,
\label{+11}
\end{equation}
and 
\begin{equation}
\xi (0)=\xi'(0)=\cdots=\xi^{(k-1)}(0)=0.
\label{+12}
\end{equation}
Then
\begin{equation}
\xi\equiv 0\quad\mbox{on}~[0,\beta].
\label{+13}
\end{equation}
\label{lemma3+}
\end{lem}
\noindent{\bf Proof.}\ We deduce, from
(\ref{+12}), that
\begin{equation}
|\xi(r)|\le Cr^{k-1+\gamma},\quad 0\le r\le \beta,
\label{+14}
\end{equation}
where $C$ is some positive constant.

Using (\ref{+14}), we deduce, from (\ref{+11}), that 
\begin{equation}
|\xi (r)|\le \frac{\alpha}{r}\int_0^r
Cs^{k-1+\gamma}~ds=\frac{C\alpha}{k+\gamma}r^{k-1+\gamma},\quad 0\le r\le\beta.
\label{+15}
\end{equation}

Using (\ref{+15}), we deduce, from (\ref{+11}), that
\[
|\xi (r)|\le \frac{\alpha}{r}\int_0^r
\frac{C\alpha}{k+\gamma}s^{k-1+\gamma}~ds=C(\frac{\alpha}{k+\gamma})^2
r^{k-1+\gamma},\quad 0\le r\le\beta.
\]

Continue this way( by induction), we have
\[
|\xi (r)|\le C(\frac{\alpha}{k+\gamma})^j r^{k-1+\gamma},\quad\forall~ 0\le
r\le\beta,~~\forall~ j=1,2,\cdots.
\]

Since $\frac{\alpha}{k+\gamma}<1$, we obtain (\ref{+13}) by sending
$j\to\infty$. Lemma~\ref{lemma3+} is established.
\vskip 5pt
\hfill $\Box$
\vskip 5pt
Continue the proof of Theorem~\ref{liouvilleR+1}. Since 
\[
1=f(\lambda_1^v(r),\cdots,\lambda_n^v(r))
=f(\lambda_1^w(r),\cdots,\lambda_n^w(r)), 
\]
we have
\begin{eqnarray*}
0&=&\int_0^1\Big(\frac{d}{dt}f(t\lambda^v(r)+(1-t)\lambda^w(r))
\Big)~dt\\
&=&\Big (\sum\limits_{i=1}^n \int_0^1
f_{\lambda_i}(t\lambda^v(r)+(1-t)\lambda^w(r))~dt\Big)
(\lambda_i^v(r)-\lambda_i^w(r)). 
\end{eqnarray*}
Since $\lambda^v(0)=\lambda^w(0)$ and
$f_{\lambda_i}(\lambda^v(0))=f_{\lambda_1}(\lambda^v(0))>0$, we deduce
from the above that
\[
\lambda_1^v(r)-\lambda^w_1(r)=-\sum\limits_{i=2}^n
(1+o(1))(\lambda_i^v(r)-\lambda^w_i(r)),
\]
where $o(1)$ denotes some quantities tending to $0$ as $r\to
0$.

Since $v'(0)=w'(0)=0$, we have 
\[
\lambda_1^v(r)-\lambda_1^w(r)=-\frac{2}{n-2}v(r)^{-\frac{n+2}{n-2}}
(v''(r)-w''(r))+O(1)(|v(r)-w(r)|+|v'(r)-w'(r)|),
\]
and, for $2\le i\le n$, 
\[
\lambda_i^v(r)-\lambda_i^w(r)=-\frac{2}{n-2}v(r)^{-\frac{n+2}{n-2}}
\frac{v'(r)-w'(r)}{r}+O(1)(|v(r)-w(r)|+|v'(r)-w'(r)|).
\]
It follows that 
\[
v''(r)-w''(r)=-\frac{n-1}{r}(v'(r)-w'(r))(1+o(1))
+O(1)(|v(r)-w(r)|+|v'(r)-w'(r)|),
\]
i.e.,
\[
(r^{n-1}(v'(r)-w'(r)))'=o(r^{n-2})|v'(r)-w'(r)|
+O(r^{n-1})(|v(r)-w(r)|+|v'(r)-w'(r)|). 
\]
Integrating the above, we have, using $v(0)-w(0)=0$,
\begin{eqnarray*}
&&|v'(r)-w'(r)|\\
&\le&\frac {o(1)}{r}\int_0^r|v'(s)-w'(s)|~ds+C\int_0^r
 (|v(s)-w(s)|+|v'(s)-w'(s)|)~ds\\
&\le& \frac {o(1)}{r}\int_0^r|v'(s)-w'(s)|~ds.
\end{eqnarray*}
Applying Lemma~\ref{lemma3+} to $\xi=v'-w'$, we have, for some $\delta>0$, 
\[
v'(r)-w'(r)\equiv 0\qquad\mbox {in}~(0,\delta).
\]
For $r\ge \delta$, the O.D.E. satisfied by $u$ and $w$ is regular, so
$v\equiv w$ in $(0,1)$. Hence $w$ is regular in $(0,1)$. Consequently,
$b\ge -1$.

\noindent
{\bf Proof of Theorem~\ref{liouvilleR+2}.}\ To give the main idea of
the proof, we first prove Theorem~\ref{liouvilleR+2} under a stronger
assumption on $u$, i.e.,
\begin{equation}
u_{0,1}(x):=|x|^{2-n}u(\frac{x}{|x|^2})~~\mbox{can be extended to a
positive function in}~~C^2(\overline{B^{+}_1}),
\label{+16a}
\end{equation}
and
\begin{equation}
A^{u_{0,1}}\in U\quad\mbox{on}~\overline{B^{+}_1}.
\label{+16b}
\end{equation}
For $x\in\Bbb {R}^n$, $\lambda>0$, let $u_{x,\lambda}$ denote the
reflection of $u$ with respect to $B_{\lambda}(x)$, i.e.,
\[
u_{x,\lambda}(y):=(\frac{\lambda}{|y-x|})^{n-2}
u(x+\frac{\lambda^2(y-x)}{|y-x|^2}).
\]
\begin{lem}
Let $u$ be as in Theorem~\ref{liouvilleR+2}. Then, for any
$x\in\partial\Bbb{R}^n_{+}$, there exists $\lambda_0(x)>0$ such that
\begin{equation}
u_{x,\lambda}\le u\quad\mbox{on}~\overline{\Bbb{R}^n_{+}}\setminus
B_{\lambda}(x),\quad \forall~0<\lambda<\lambda_0(x).
\label{+17}
\end{equation}
\label{lemma2.1+}
\end{lem}
\noindent{\bf Proof of Lemma~\ref{lemma2.1+}.}\ We follow the arguments
in the proof of lemma~2.1 in \cite {LZ}. Without loss of generality, take
$x=0$ in (\ref{+17}), and we use $u_{\lambda}$ to denote
$u_{0,\lambda}$. By the $C^1$ regularity of $u$, there exists $r_0>0$
such that
\[
\frac{d}{dr}(r^{\frac{n-2}{2}}u(r,\theta))>0,\quad\forall~0<r<r_0,~\theta\in
S^{n-1},
\] 
from which, we deduce
\begin{equation}
u_{\lambda}(y)<u(y),\quad\forall~0<\lambda<|y|<r_0.
\label{+18}
\end{equation}
Because of (\ref{+6b}), there exists some constant $\alpha>0$ such that
\[
u(y)\ge \frac{\alpha}{|y|^{n-2}},\quad\forall~|y|\ge r_0.
\]
Let
$\lambda_0=\min\{\alpha^{\frac{1}{n-2}}(\max\limits_{\overline{B^{+}_{r_0}}
}u)^{\frac{1}{2-n}},r_0\}$. Then
\[
u_{\lambda}(y)\le
(\frac{\lambda_0}{|y|})^{n-2}(\max\limits_{\overline{B^{+}_{r_0}}}u)\le
\frac{\alpha}{|y|^{n-2}}\le
u(y),\quad\forall~0<\lambda<\lambda_0,~|y|\ge r_0.
\]
(\ref{+17}) with $x=0$ follows from (\ref{+18}) and the
above. Lemma~\ref{lemma2.1+} is established.
\vskip 5pt
\hfill $\Box$
\vskip 5pt
For $x\in\partial\Bbb{R}^n_{+}$, let 
\begin{equation}
\bar\lambda(x):=\sup\{\mu>0|u_{x,\lambda}\le
u~~\mbox{on}~~\overline{\Bbb{R}^n_{+}}\setminus
B_{\lambda}(x),~~~\forall~0<\lambda<\mu\}. 
\label{+9-1}
\end{equation}
Clearly, $\bar\lambda(x)>0$. On the other hand,
$\bar\lambda(x)<\infty$ because of (\ref{+6b}).
\begin{lem}
Let $u$ be as in Theorem~\ref{liouvilleR+2}, and we further assume that
$u$ satisfies (\ref{+16a}) and (\ref{+16b}). Then, for all
$x\in\partial\Bbb{R}^n_{+}$,  
\begin{equation}
u_{x,\bar\lambda(x)}\equiv u\quad\mbox{on}~\Bbb{R}^n_{+}\setminus\{x\}.
\label{+19}
\end{equation}
\label{lemma2.2+}
\end{lem}
\noindent{\bf Proof of Lemma~\ref{lemma2.2+}.}\ Without loss of
generality, take $x=0$. We use notation $\bar\lambda=\bar\lambda(0)$
and $u_{\lambda}=u_{0,\lambda}$. By the definition of $\bar\lambda$,
\begin{equation}
u_{\bar\lambda}\le u\quad\mbox{on}~\overline{\Bbb{R}^n_{+}}\setminus
B_{\bar\lambda}.
\label{+20}
\end{equation}
From now on, we always assume that (\ref{+19}) does not hold for $x=0$,
and we will reach a contradiction. We first show that
\begin{equation}
u-u_{\bar\lambda}>0\quad\mbox{on}~
\overline{\Bbb{R}^n_{+}}\setminus\overline{B^{+}_{\bar\lambda}}.
\label{+21}
\end{equation}
Indeed, if, for some $\bar
x\in\Bbb{R}^n_{+}\setminus\overline{B^{+}_{\bar\lambda}}$,
$(u-u_{\bar\lambda})(\bar x)=0$. Using (\ref{+6a}) and
hypotheses (\ref{U1}) and (\ref{F1}), we have
\[
F(A^{u_{\bar\lambda}})=1
\quad\mbox{on}~\Bbb{R}^n_{+}\setminus\overline{B_{\bar\lambda}}.
\]
A calculation yields, using (\ref{+6a}),
\[
\frac{\partial u_{\bar\lambda}}{\partial
x_n}=cu_{\bar\lambda}^{\frac{n}{n-2}},\quad
\mbox{on}~\partial\Bbb{R}^n_{+}\setminus B_{\bar\lambda}.
\]
Arguing as in the proof of lemma~2.1 in \cite{LL} (using hypotheses
(\ref{U2}) and (\ref{F2})), we have, near $\bar x$,
\begin{equation}
0=F(A^u)-F(A^{u_{\bar\lambda}})=L(u-u_{\bar\lambda}),
\label{+22}
\end{equation}
where $L=-a_{ij}(x)\partial_{ij}+b_i(x)\partial_i+c(x)$ is an elliptic
operator with continuous coefficients. By the strong maximum
principle, $u-u_{\bar\lambda}\equiv 0$ near $\bar x$. This implies
(\ref{+19}) for $x=0$, a contradiction.

If $(u-u_{\bar\lambda})(\bar x)=0$ for some
$\bar x\in\partial\Bbb{R}^n_{+}\setminus\overline{B^{+}_{\bar\lambda}}$,
we have
\[
\frac{\partial(u-u_{\bar\lambda})}{\partial x_n}(\bar
x)=(cu^{\frac{n}{n-2}}-cu_{\bar\lambda}^{\frac{n}{n-2}})(\bar x)=0.
\]
Since we still have (\ref{+22}) near $\bar x$, we apply the Hopf Lemma
to obtain that $u-u_{\bar \lambda}\equiv 0$ near $\bar x$, again
leading to (\ref{+19}) for $x=0$, a contradiction. We have established
(\ref{+21}). 

Next we show that 
\begin{equation}
\lim_{y\in \bar R^n_+, |y|\to\infty}
|y|^{n-2}(u(y)-u_{\bar \lambda }(y))>0.
\label{+23}
\end{equation}
Let $x=\frac{y}{|y|^2}$, we have
\[
|y|^{n-2}u(y)=u_{0,1}(x),\quad
|y|^{n-2}u_{\bar\lambda}(y)=\bar\lambda^{n-2}u(\frac{\bar\lambda^2y}{|y|^2})
=\bar\lambda^{n-2}u(\bar\lambda^2x)=:v(x).
\]
By (\ref{+16a}), (\ref{+16b}) and the conformal invariance of the
equation (\ref{+6a}), both $u_{0,1}$ and $v$ are $C^2$ solutions of
(\ref{+6a}). We also know, from (\ref{+21}), that 
\[
u_{0,1}-v>0\quad\mbox{in}~B^{+}_{\frac{1}{\bar\lambda}}.
\]
By the same arguments used in proving $u-u_{\bar\lambda}>0$ on
$\partial\Bbb{R}^n_{+}\setminus\overline{B^{+}_{\bar\lambda}}$, we
have 
\[
(u_{0,1}-v)(0)>0,
\]
which implies (\ref{+23}).

Since $u-u_{\bar\lambda}=0$ on $\partial
B_{\bar\lambda}\cap\Bbb{R}^n_{+}$ and (\ref{+21}), we can apply the
Hopf Lemma as in the proof of lemma~2.1 in \cite {LL} (see also the outlines
near (\ref{+22})) to obtain
\begin{equation}
\frac{\partial
(u-u_{\bar\lambda})}{\partial\nu}>0\quad\mbox{on}~\partial
B_{\bar\lambda}\cap\Bbb{R}^n_{+}, 
\label{+24}
\end{equation}  
where $\nu$ denotes the unit outer normal of $\partial
B_{\bar\lambda}$.

At last we prove that 
\begin{equation}
\frac{\partial
(u-u_{\bar\lambda})}{\partial\nu}>0\quad\mbox{on}~\partial
B_{\bar\lambda}\cap\partial \Bbb{R}^n_{+}, 
\label{+25}
\end{equation}
where $\nu$ still denotes the unit outer normal of $\partial
B_{\bar\lambda}$. 

Let $\bar x\in\partial B_{\bar\lambda}\cap\partial\Bbb{R}^n_+$. Then
as in the proof of lemma~2.1 in \cite {LL}, we have (\ref{+22}) near
$\bar x$ with continuous coefficients. Clearly, for some constant $A>0$,  
\[
|\frac{\partial (u-u_{\bar\lambda})}{\partial
x_n}|=|c(u^{\frac{n}{n-2}}-u_{\bar\lambda}^{\frac{n}{n-2}})|\le A
(u-u_{\bar\lambda}),
\quad\mbox{in}~(\Bbb{R}^n_{+}\setminus\overline{B^{+}_{\bar\lambda}})\cap
B_1(\bar x).
\]
By (\ref{+22}), and for a possibly larger $A$, we have
\[
a_{ij}\partial_{ij}(u-u_{\bar\lambda})+b_i\partial_i(u-u_{\bar\lambda})\le
A(u-u_{\bar\lambda}),\quad\mbox{in}~
(\Bbb{R}^n_{+}\setminus\overline{B^{+}_{\bar\lambda}})\cap B_1(\bar x).
\]
Now an application of lemma~10.1 in \cite {LZ} (with $\Omega
=(\Bbb{R}^n_{+}\setminus\overline{B_{\bar\lambda}})\cap B_1(\bar
x)$, $\sigma =x_n$, $\rho=|x|^2-\bar\lambda^2$, and our
$u-u_{\bar\lambda}$ being the $u$ there) yields 
\[
\frac{\partial(u-u_{\bar\lambda})}{\partial \nu}(\bar x)>0.
\]
So we have established (\ref{+25}).

Given (\ref{+21}), (\ref{+23}), (\ref{+24}), (\ref{+25}), the
positivity and continuity of $u$ on $\overline{\Bbb{R}^n_{+}}$, we can
easily prove that there exists some $\epsilon>0$ such that
\[
u_{\lambda}\le u\quad\mbox{on}~\overline{\Bbb{R}^n_{+}}\setminus
B^{+}_{\lambda},\quad \forall~\bar\lambda\le\lambda\le\bar\lambda+\epsilon,
\]
which violates the definition of $\bar\lambda$. Lemma~\ref{lemma2.2+} is
established. 
\vskip 5pt
\hfill $\Box$
\vskip 5pt

{\bf The Proof of Theorem~\ref{liouvilleR+2} under the additional
hypotheses (\ref{+16a}) and (\ref{+16b}).}\ Let $u$ be as in
Theorem~\ref{liouvilleR+2} and $u$ satisfies (\ref{+16a}) and
(\ref{+16b}). By Lemma~\ref{lemma2.2+} and a calculus lemma used in
\cite {LZhu} (see, e.g., lemma~11.1 in \cite {LZ}), 
\begin{equation}
u(x',0)\equiv \frac{\hat a}{(|x'-\bar x'|^2+d^2)^{\frac{n-2}{2}}},
\quad\mbox{on}~\Bbb{R}^{n-1},
\label{+26}
\end{equation}
where $\bar x'\in\Bbb{R}^{n-1}$, $\hat a$ and $d$ are positive
constants.

Let $P=(\bar x',-d)$ and define
\[
v(z):=(\frac{2d}{|z-P|})^{n-2}u(P+\frac{4d^2(z-P)}{|z-P|^2}).
\]
By the arguments in \cite {LZhu} and \cite {Bi}, as in the proof of
lemma~4.5 in \cite {LZ}, we know that $v$ is radially symmetric with respect
to $Q:=(\bar x',d)$ in $B_{2d}(Q)$. By the conformal invariance of the
equation satisfied by $u$, we have 
\[
F(A^v)=1,\quad A^v\in U,\quad v>0,\quad\mbox{in}~\overline{B_{2d}(Q)}.
\]
By Theorem~\ref{liouvilleR+1},
\[
v(z)\equiv (\frac{\bar a}{1+\bar
b|z-Q|^2}^{\frac{n-2}{2}}\quad\mbox{in}~\overline{B_{2d}(Q)}, 
\]
where $\bar a>0$ and $1+\bar b(2d)^2>0$. Compare this to (\ref{+26}),
we must have $\bar b>0$. This, together with (\ref{+26}), implies 
\[
u(x)\equiv (\frac{a}{1+b|x-\bar
x|^2})^{\frac{n-2}{2}}\quad\mbox{on}~\Bbb{R}^n_{+}, 
\]
where $a=d^{-2}\hat a^{\frac{2}{n-2}}$, $b=d^{-2}$, $\bar x=(\bar
x',\bar x_n)$, $\hat a,~d,~\bar x'$ are given in (\ref{+26}), and $\bar
x_n$ is some real number.

Since $A^u(0)=2a^{-2}bI$, we have $2a^{-2}bI\in U$ and
$F(2a^{-2}bI)=F(A^u(0))=1$. By the boundary condition of $u$ at
$x=0$, we have $(n-2)a^{-1}b\bar x_n=c$. Theorem~\ref{liouvilleR+2} is
established under the additional hypotheses. 
\vskip 5pt
\hfill $\Box$
\vskip 5pt

{\bf The proof of Theorem~\ref{liouvilleR+2}.}\ By Lemma~\ref{lemma2.1+},
there exists $\lambda_0>0$ such that
\begin{equation}
u_{\lambda}\le u\quad\mbox{on}~\overline{\Bbb{R}^n_{+}}\setminus
B_{\lambda},\quad \forall~0<\lambda<\lambda_0,
\label{+1-1}
\end{equation}
where $u_{\lambda}=u_{0,\lambda}$ and $B_{\lambda}=B_{\lambda}(0)$.

Let $w=u_{0,1}$. As in the proof of Lemma~\ref{lemma2.1+} and in the
proof of lemma~2.1 of \cite {LL}, there exists some $\lambda_1>0$ such that 
\begin{equation}
w_{\lambda}\le w\quad\mbox{on}~\overline{\Bbb{R}^n_{+}}\setminus
B_{\lambda},\quad \forall~0<\lambda<\lambda_1.
\label{+2-1}
\end{equation}
Rewriting (\ref{+1-1}) and (\ref{+2-1}) as 
\[
w_{\lambda}\le
w\quad\mbox{in}~B_{\lambda}^+,\quad\forall~\lambda>\frac{1}{\lambda_0} ,
\]
and
\[
w_{\lambda}\ge w\quad\mbox{in}~B_{\lambda}^+,\quad\forall~0<\lambda<\lambda_1.
\]
Let
\[
\underline{\lambda}:=\sup\{\mu|~w_{\lambda}(x)\ge
w(x),\quad\forall~0<|x|\le\lambda\le \mu\},
\]
and
\[
\overline{\lambda}:=\inf\{\mu|~w_{\lambda}(x)\le
w(x),\quad\forall~\lambda\ge\mu,~0<|x|\le\lambda\}.
\]
If $\overline{\lambda}\le\underline{\lambda}$, then
$w_{\underline{\lambda}}\equiv w_{\overline{\lambda}}\equiv w$, and
$u$ satisfies (\ref{+16a}) and (\ref{+16b}). Theorem~\ref{liouvilleR+2} in
this case has already been established. In the following, we assume
that $\overline{\lambda}>\underline{\lambda}$ and we will reach a
contradiction. \newline
Clearly, $w_{\lambda}(0)=\frac{1}{\lambda^{n-2}}u(0)$, so we have
\[
\frac{1}{\overline{\lambda}^{n-2}}u(0)\le
w(0)\le\frac{1}{\underline{\lambda}^{n-2}}u(0). 
\]
Since $\overline{\lambda}>\underline{\lambda}$, there must be at least
one strict inequality in the above. Without loss of generality, we
assume that 
\[
w_{\overline{\lambda}}(0)=\frac{1}{\overline{\lambda}^{n-2}}u(0)<w(0).
\]
This guarantees that there is no touching of $w_{\lambda}$ and $w$
near $0$ for $\lambda$ close to $\overline{\lambda}$. Therefore, by
the moving sphere arguments used earlier, we have, for $\lambda$ close
to $\overline{\lambda}$, that 
\[
w_{\lambda}\le w\quad\mbox{in}~B_{\lambda}.
\]
This violates the definition of
$\overline{\lambda}$. Theorem~\ref{liouvilleR+2} is established.
\vskip 5pt
\hfill $\Box$
\vskip 5pt
\subsection{Proof of Theorem~\ref{liouvilleR+3}}
\noindent{\bf Proof of Theorem~\ref{liouvilleR+3}.}\ Let 
\[ 
\alpha:=\liminf\limits_{x\to\Bbb{R}^n_{+},|x|\to\infty}(|x|^{n-2}u(x))\in 
[0,\infty]. 
\] 
\begin{lem} 
\[ 
\alpha>0. 
\] 
\label{lemma4-1+} 
\end{lem} 
{\bf Proof of Lemma~\ref{lemma4-1+}.}\ We follow the arguments of the 
proof of lemma~4.1 in \cite {LZ}. Let 
\[ 
O:=\{y\in\Bbb{R}^n_{+}|u(y)<|y|^{2-n}\}. 
\] 
To prove the lemma, we only need to show  
\[ 
\liminf\limits_{x\in O,|x|\to\infty}|x|^{n-2}u(x)>0. 
\] 
We know 
\[ 
\Delta u\le 0\quad\mbox{in}~O, 
\] 
and  
\[ 
\frac{\partial u}{\partial x_n}=cu^{\frac{n}{n-2}}\le
(|c|+1)|y|^{-2}u\quad\mbox{on}~\partial O\cap\partial\Bbb{R}^n_{+}.  
\] 
For $A>1$, let  
\[ 
\xi (y):=|y-Ae_n|^{2-n}+|y|^{1-n}. 
\] 
For large $A$ and $R=A^2$, we have 
$$ 
\left\{ 
\begin{array}{lcl} 
-\Delta\xi&\le&0,\quad\mbox{on}~\Bbb{R}^n_{+}\setminus B_R\\ 
\frac{\partial\xi}{\partial x_n}(y)&\ge&\frac{|c|+1}{|y|^2}\xi(y),\quad 
y\in\partial\Bbb{R}^n_{+}\setminus B_R. 
\end{array} 
\right. 
$$
Take $\bar\epsilon(A)>0$ be a small constant such that 
\[ 
w:=u-\bar\epsilon\xi\ge 0\quad\mbox{on}~\partial(O\setminus 
B_R)\cap\Bbb{R}^n_{+}.  
\] 
If follows that 
$$
\left\{ 
\begin{array}{lcl} 
\Delta w&\le&0,\quad\mbox{on}~O\setminus B_R\\ 
\frac{\partial w}{\partial x_n}(y)-\frac{|c|+1}{|y|^2}w(y)&\le&0,\quad 
\forall~y\in\partial (O\setminus B_R)\cap\partial\Bbb{R}^n_{+}. 
\end{array} 
\right. 
$$
Clearly, $\liminf\limits_{x\in O\setminus B_R,|x|\to \infty}w(x)\ge 
0$. By the maximum principle,  
\[ 
w\ge 0\quad\mbox{on}~O\setminus B_R. 
\] 
Hence 
\[ 
\liminf\limits_{x\in O,|x|\to\infty}|x|^{n-2}u(x)\ge\bar\epsilon>0. 
\] 
Lemma~\ref{lemma4-1+} is proved. 
\vskip 3pt 
\hfill $\Box$ 
\vskip 3pt 
\begin{lem} 
For any $x\in\partial\Bbb{R}^n_{+}$, there exists $\lambda_0(x)>0$ 
such that 
\[ 
u_{x,\lambda}\le u\quad\mbox{on}~\Bbb{R}^n_{+}\setminus 
B_{\lambda}(x),\quad \forall~0<\lambda<\lambda_0(x). 
\] 
\label{lemma4-2+} 
\end{lem} 
{\bf Proof of Lemma~\ref{lemma4-2+}.}\ Since we know 
$\alpha>0$. Lemma~\ref{lemma4-2+} follows from the proof of 
Lemma~\ref{lemma2.1+}. 
\vskip 3pt 
\hfill $\Box$ 
\vskip 3pt 
For $x\in\partial\Bbb{R}^n_{+}$, let $\bar\lambda (x)$ be defined as 
in (\ref{+9-1}). By Lemma~\ref{lemma4-2+}, $\bar\lambda (x)>0$. 
\begin{lem} 
If $\alpha=\infty$, then  
\[ 
\bar\lambda (x)=\infty,\quad\forall~x\in\partial\Bbb{R}^n_{+}. 
\] 
If $\alpha<\infty$, then 
\begin{equation} 
\bar\lambda (x)^{n-2}u(x)=\alpha,\quad\forall~x\in\partial\Bbb{R}^n_{+}. 
\label{more3-1} 
\end{equation} 
\label{lemma4-3+} 
\end{lem} 
{\bf Proof of Lemma~\ref{lemma4-3+}.}\ By the definition of $\bar 
\lambda(x)$, we know 
\[ 
u_{x,\lambda}(y)\le u(y),\quad\forall~0<\lambda<\bar\lambda(x), 
~~~~\forall~y\in\Bbb{R}^n_{+}\setminus B_{\lambda}(x). 
\] 
It follows that  
\[ 
\lambda^{n-2}u(x)=\liminf\limits_{y\in\Bbb{R}^n_{+},|y|\to\infty} 
|y|^{n-2}u_{x,\lambda}(y)\le\liminf\limits_{y\in\Bbb{R}^n_{+},|y|\to\infty} 
|y|^{n-2}u(y)=\alpha,\quad\forall~0<\lambda<\bar\lambda(x), 
\] 
If $\alpha<\infty$, we have 
\[ 
\bar\lambda(x)^{n-2}u(x)\le\alpha<\infty, 
\quad\forall~x\in\partial\Bbb{R}^n_{+}. 
\] 
In fact we must have  
\[ 
\bar\lambda (x)^{n-2} u(x)=\alpha. 
\] 
Indeed, if $\bar\lambda (x)^{n-2}u(x)<\alpha$, then  
\[ 
\lim\limits_{y\in\Bbb{R}^n_{+},|y|\to\infty}|y|^{n-2}(u(y)-u_{x,\bar\lambda 
(x)}(y))=\alpha-\bar\lambda (x)^{n-2} u(x)>0, 
\] 
and the arguments in the proof of Lemma~\ref{lemma2.2+} show that the 
moving sphere procedure should not stop at $\bar\lambda (x)$, 
violating the definition of $\bar\lambda (x)$. 
 
Now assume $\alpha=\infty$. Without loss of generality, we show 
$\bar\lambda:=\bar\lambda(0)=\infty$. We prove it by contradiction. Suppose 
$\bar\lambda<\infty$. By the definition of $\bar\lambda$, (\ref{+20}) 
holds. Since $\alpha=\infty$, we have 
\[ 
\liminf\limits_{y\in\overline{\Bbb{R}^n_{+}},|y|\to\infty} 
(u(y)-u_{\bar\lambda}(y))|y|^{n-2}=\infty. 
\] 
This plays the same role as (\ref{+23}) in the proof of  
Lemma~\ref{lemma2.2+}, and the arguments there lead to a contradiction 
to the definition of $\bar \lambda$. Lemma~\ref{lemma4-3+} is established. 
\vskip 3pt 
\hfill $\Box$ 
\vskip 3pt  
To prove Theorem~\ref{liouvilleR+3}, we first consider the case 
$\alpha<\infty$. Our proof goes along the line of the proof of 
Theorem~\ref{theorem1}. Our next lemma, whose proof is given towards
the end of this section, is an analogue of Lemma~\ref{lemma0}.  
\begin{lem} 
For $n\ge 3$, $a,d>0$, $c\in\Bbb{R}$, $p,q\in\Bbb{R}^{n-1}$ and $p\neq 
q$, let $u\in  
C^1(\overline{B_d^{+}}\setminus\{0\})$ satisfy 
\begin{equation} 
\left\{ 
\begin{array}{lcl} 
\Delta u\le 0,\quad~~~~~~\mbox{in}~B_d^{+}~\mbox{in 
the distribution sense},&&\\ 
\frac{\partial u}{\partial 
  x_n}=cu^{\frac{n}{n-2}},\quad\mbox{on}~(\partial 
B_d^{+}\cap\partial\Bbb{R}^n_{+})\setminus\{0\} ,&&\\ 
u(x)\ge\max\{a+p\cdot 
x'+ca^{\frac{n}{n-2}}x_n-\bar\delta(|x|),&&\\
~~~~~~~~~~~~~~~~a+q\cdot
x'+ca^{\frac{n}{n-2}}x_n-\bar\delta(|x|)\},\quad\forall~x\in
B_d^{+},&&   
\end{array} 
\right. 
\label{A4e1+} 
\end{equation} 
where $x'=(x_1,\cdots,x_{n-1})$, $\bar\delta (r)>0$ and 
$\lim\limits_{r\to 0^{+}}\frac{\bar\delta (r)}{r}=0$. 
Then  
\[ 
\liminf\limits_{x\in B_1^{+},x\to 0} u(x)> a. 
\] 
\label{lemA4+} 
\end{lem} 

\begin{lem} 
Under the hypothesis of Theorem~\ref{liouvilleR+3}, if $\alpha<\infty$, 
then $u$ is of the form (\ref{11new}) with $\bar x$, $a$ and $b$ given 
below (\ref{11new}).  
\label{proposition4-1+} 
\end{lem} 
{\bf Proof of Lemma~\ref{proposition4-1+}.}\ For 
$x\in\partial\Bbb{R}^n_{+}$, let 
\[ 
\phi^{(x)}(y):=x+\frac{\bar\lambda (x)^2(y-x)}{|y-x|^2},\quad \psi 
(y)=\frac{y}{|y|^2},\quad 
w^{(x)}(y):=(u_{\phi^{(x)}})_{\psi}=u_{\phi^{(x)}\circ\psi}. 
\] 
By the definition of $\bar\lambda (x)$, 
\begin{equation} 
u\ge u_{\phi^{(x)}}\quad\mbox{on}~\Bbb{R}^n_{+}\setminus 
B_{\bar\lambda (x)}(x),\quad\forall~x\in\partial\Bbb{R}^n_{+}. 
\label{17new} 
\end{equation} 
By (\ref{more3-1}), 
\[ 
w^{(x)}(0):=\bar\lambda(x)^{n-2}u(x)=\alpha, 
\quad\forall~x\in\partial\Bbb{R}^n_{+}. 
\]  
We have 
$$
\left\{ 
\begin{array}{lcl} 
u_{\psi}\in C^2(\Bbb{R}^n_{+}),&&\\ 
\Delta u_{\psi}\le 0\quad\mbox{in}~\Bbb{R}^n_{+}~~\mbox{since}~~\Delta 
u\le 0~~\mbox{in}~~\Bbb{R}^n_{+},&&\\ 
\liminf\limits_{\Bbb{R}^n_{+}\ni y\to 
0}u_{\psi}(y)=\liminf\limits_{z\in\Bbb{R}^n_{+},|z|\to\infty}|z|^{n-2}u(z) 
=\alpha,&& 
\end{array} 
\right. 
$$
and it is clear, for some $\delta (x)>0$ and by (\ref{17new}), that 
\begin{eqnarray*} 
w^{(x)}\in 
C^2(\overline{B^{+}_{\delta(x)}}),&&\quad\forall~x\in\partial\Bbb{R}^n_{+},\\ 
u_{\psi}\ge w^{(x)},&&\quad\mbox{in}~B^{+}_{\delta (x)}. 
\end{eqnarray*} 
By (\ref{+6a}) and the conformal invariance of the boundary condition 
satisfied by $u$, 
$$
\left\{ 
\begin{array}{lcl} 
\frac{\partial w^{(x)}}{\partial y_n}&=&c[w^{(x)}]^{\frac{n}{n-2}}, 
\quad\mbox{on}~\partial\Bbb{R}^n_{+}\setminus\{\frac{x}{|x|^2}\} \\ 
\frac{\partial u_{\psi}}{\partial y_n}&=&c[u_{\psi}]^{\frac{n}{n-2}}, 
\quad\mbox{on}~\partial\Bbb{R}^n_{+}\setminus\{0\}. 
\end{array} 
\right. 
$$
By Lemma~\ref{lemA4+},  
\[ 
\nabla_{y'}w^{(x)}(0)=\nabla_{y'}w^{(0)}(0), 
\quad\forall~x\in\partial\Bbb{R}^n_{+}. 
\] 
So for $x=(x',0)$, 
\begin{eqnarray*} 
\vec{V}:&=&\nabla_{y'}w^{(0)}(0)=(n-2)\bar\lambda(x)^{n-2}u(x)x 
+\bar\lambda(x)^n\nabla_{x'}u(x)+\bar\lambda(x)^n\nabla_{x'}u(x)\\ 
&=&(n-2)\alpha x'+\alpha^{\frac{n}{n-2}}u(x)^{\frac{n}{2-n}}\nabla_{x'}u(x). 
\end{eqnarray*} 
Thus we have 
\[ 
\nabla_{x'}[\frac{n-2}{2}\alpha^{\frac{n}{n-2}}u(x',0)^{-\frac{2}{n-2}} 
-\frac{n-2}{2}|x'|^2+\vec{V}\cdot x']=0, 
\] 
which implies, for some $\bar x'\in\Bbb{R}^{n-1}$, and $d\in R$, that 
\[ 
u(x',0)^{-\frac{2}{n-2}}\equiv \alpha^{-\frac{2}{n-2}}|x'-\bar 
x'|^2+d\alpha^{-\frac{2}{n-2}}. 
\] 
Since $u>0$, we have $d>0$ and 
\begin{equation} 
u(x',0)\equiv \Big(\frac{\alpha^{\frac{2}{n-2}}}{d+|x'-\bar 
x'|^2}\Big)^{\frac{n-2}{2}}.  
\label{18new} 
\end{equation} 
For simplicity, we take $\bar x'=0$. By (\ref{more3-1}) and the above, 
\[ 
\alpha=\bar\lambda(0)^{n-2}u(0)=\bar\lambda 
(0)^{n-2}\frac{\alpha}{d^{\frac{n-2}{2}}}, 
\] 
which gives $\bar\lambda :=\bar\lambda (0)=\sqrt{d}$.\newline 
Since  
\[ 
u_{\bar\lambda}(y)=(\frac{\bar\lambda}{|y|})^{n-2}u(\frac{\bar\lambda^2 
y}{|y|^2}), 
\] 
we have, by (\ref{18new}), 
\[ 
u_{\bar\lambda}(x',0)=\frac{\bar\lambda^{n-2}\alpha} 
{(d|x'|^2+\bar\lambda^4)^{\frac{n-2}{2}}} 
=\frac{\alpha}{(|x'|^2+d)^{\frac{n-2}{2}}} 
=u(x',0),\quad\forall~x'\in\Bbb{R}^{n-1}. 
\] 
Thus by the conformal invariance of the equation and the boundary 
condition satisfied by $u$, we have 
$$
\left\{ 
\begin{array}{lcl} 
F(A^u)=F(A^{u_{\bar\lambda}})=1,~~A^u\in U,~~A^{u_{\bar\lambda}}\in
U,\quad&&\mbox{in}~\overline {\Bbb{R}^n_{+}}\setminus \{0\},\\ 
u-u_{\bar\lambda}=0,\quad&&\mbox{on}~\partial\Bbb{R}^n_{+}\setminus \{0\},\\ 
\frac{\partial (u-u_{\bar\lambda})}{\partial 
x_n}=cu^{\frac{n}{n-2}}-cu_{\bar\lambda}^{\frac{n}{n-2}}=0, 
\quad&&\mbox{on}~\partial\Bbb{R}^n_{+}\setminus \{0\},\\  
u-u_{\bar\lambda}\ge 0,\quad&&\mbox{on}~\Bbb{R}^n_{+}\setminus B_{\bar\lambda}.
\end{array} 
\right. 
$$
As usual, $u-u_{\bar\lambda}$ satisfies a linear second order elliptic 
equation and therefore, by the Hopf lemma and the strong maximum principle,  
\[ 
u-u_{\bar\lambda}\equiv 0\quad\mbox{on}~\Bbb{R}^n_{+}. 
\] 
In particular, $u$ satisfies (\ref{+16a}) and (\ref{+16b}). So $u$ is of the 
form (\ref{11new}) by our earlier discussion of Theorem~\ref{liouvilleR+3} 
under (\ref{+16a}) and (\ref{+16b}). Lemma~\ref{proposition4-1+} is 
established.  
\vskip 3pt 
\hfill $\Box$ 
\vskip 3pt 
\begin{lem} 
Under the hypothesis of Theorem~\ref{liouvilleR+3} except (\ref{more1}), 
if $\alpha=\infty$, then   
\begin{equation} 
u(x',x_n)\equiv 
u(0',x_n),\quad\forall~x'\in\Bbb{R}^{n-1},~\forall~x_n\ge 0. 
\label{19new} 
\end{equation} 
Moreover $c\ge 0$, and if $c=0$, $u$ must be a constant. 
\label{prop2new} 
\end{lem} 
\noindent{\bf Proof of Lemma~\ref{prop2new}.}\ Since $\alpha=\infty$, we 
have, by Lemma~\ref{lemma4-3+}, 
\[ 
\bar\lambda(x)=\infty,\quad\forall~x\in\partial\Bbb{R}^n_{+}. 
\] 
i.e., 
\[ 
u_{x,\lambda}\le u \quad\mbox{on}~\Bbb{R}^n_{+}\setminus 
B_{\lambda}(x),\quad\forall~0<\lambda<\infty, 
\] 
which, by a calculus lemma(see, e.g., lemma~11.3 in \cite {LZ}), implies 
(\ref{19new}). Let  
\[ 
h(t):=u(0',t)\quad\mbox{for}~ t\ge 0. 
\]  
Since $\Delta u\le 0$, we have 
\[ 
h''(t)\le 0,\quad\forall~t\ge 0, 
\]  
so  
\[ 
h'(t)\le h'(s),\quad\forall~t\ge s\ge 0. 
\] 
Hence 
\[ 
h(t)-h(s)\le h'(s)(t-s),\quad\forall~t\ge s\ge 0.  
\] 
and 
\[ 
h'(s)\ge \liminf\limits_{t\to\infty}\frac{h(t)-h(s)}{t-s}\ge 
0,\quad\forall~s\ge 0.  
\]  
Since $\frac{\partial u}{\partial x_n}=cu^{\frac{n}{n-2}}$ on 
$\partial\Bbb{R}^n_{+}$,  
\[ 
h'(0)=ch(0)^{\frac{n}{n-2}}. 
\] 
Since $h(0)>0$ and $h'(0)\ge 0$, we have $c\ge 0$. If $c=0$, we have 
$h'(0)=0$. Recall that $h''(t)\le 0$, so  
\[ 
h'(t)\le h'(0)=0,\quad\forall~t\ge 0. 
\] 
On the other hand, $h'(t)\ge 0$, so $h'(t)\equiv 0$ and $h(t)\equiv 
h(0)$. Lemma~\ref{prop2new} is established.\newline  
{\bf Proof of Theorem~\ref{liouvilleR+3}.}\ If $\alpha<\infty$, the
theorem follows from Lemma~\ref{proposition4-1+}. If $\alpha =\infty$, then by
Lemma~\ref{prop2new}, (\ref{19new}) holds, and we only need to  
rule out the possibility of $c>0$. For this aim, we make use of 
(\ref{more1}). As before, let 
\[ 
h(t):=u(0',t),\quad\forall~t\ge 0.  
\] 
{\bf Claim.}\  $\forall~a>0$, 
\begin{equation} 
\lim\limits_{t\to\infty}\frac{h'(t)}{h(t)^a}=0. 
\label{+55} 
\end{equation} 
Indeed, if $\lim\limits_{t\to\infty}h(t)=\infty$, then (\ref{+55}) is 
obvious by $0\le h'(t)\le h'(0)$. Otherwise, there exists some 
$b\in [h(0),\infty)$ such that 
\[ 
\lim\limits_{t\to\infty}h(t)=b. 
\] 
We also know that $\lim\limits_{t\to\infty}h'(t)$ exists since 
$h''(t)\le 0$. So, by the boundedness of $h(t)$, we must have  
\[ 
\lim\limits_{t\to\infty}h'(t)=0, 
\] 
which yields (\ref{+55}). 
 
Let $(\lambda_1,\lambda_2,\cdots,\lambda_n)$ denote the eigenvalues of 
$A^u$. Then  
\[ 
\left\{ 
\begin{array}{clc} 
\lambda_1(t)=\cdots=\lambda_{n-1}(t) 
=-\frac{2}{(n-2)^2}\frac{h'(t)^2}{h(t)^{\frac{2n}{n-2}}},&&\\ 
\lambda_n(t)=-\frac{2}{n-2}\frac{h''(t)}{h(t)^{\frac{n+2}{n-2}}}+\frac{2(n-1)} 
{(n-2)^2}\frac{h'(t)}{h(t)^{\frac{2n}{n-2}}}. 
\end{array} 
\right. 
\] 
By (\ref{+55}) and the equation satisfied by $u$, 
$$ 
\left\{ 
\begin{array}{lcl} 
f(\lambda_1,\lambda_2,\cdots,\lambda_n)=1,&&\\ 
\lambda_1=o(1),\cdots,\lambda_{n-1}=o(1),&&\\ 
\lambda_n=-\frac{2}{n-2}h^{-\frac{n+2}{n-2}}h''+o(1). 
\end{array} 
\right. 
$$
By assumption (\ref{more1}), there exists some $\delta>0$ such that 
\[ 
|(\lambda_1,\cdots,\lambda_n)|\ge\delta, 
\] 
so for large $t$, 
\[ 
-\frac{2}{n-2}\frac{h''(t)}{h(t)^{\frac{n+2}{n-2}}}\ge\frac{\delta}{2}, 
\] 
i.e., 
\[ 
-h''(t)\ge\frac{n-2}{4}\delta 
h(t)^{\frac{n+2}{n-2}}>\frac{n-2}{4}\delta h(0)^{\frac{n+2}{n-2}}.   
\] 
Integrating the above inequality twice leads to  
\[ 
-h(t)+h(0)+h'(0)t\ge\frac{n-2}{8}\delta 
h(0)^{\frac{n+2}{n-2}}t^2,\quad\forall~t\ge 0. 
\] 
Sending $t\to\infty$ in the above yields a contradiction to the 
positivity of $h$. Thus we have ruled out the possibility that
$c>0$. Theorem~\ref{liouvilleR+3} is established.  
\vskip 3pt 
\hfill $\Box$ 
\vskip 3pt  
In the rest of this section, we prove Lemma~\ref{lemA4+}. We use 
notations 
\[
e_1=(1,0\cdots,0),~x=(x_1,\cdots,x_n)=(x',x_n),~B_r^{+}
=B_r\cap\Bbb{R}^n_{+},~\mbox{and}~\partial' B_1^+:=\partial
B_1^+\cap\Bbb {R}^n_+. 
\] 
Fixing some small $b>0$ to be specified later, let  
\[ 
\phi_b(x)= 
\left\{ 
\begin{array}{lcl} 
x_1,\quad\forall~x\in\partial 
B_1\cap\{x|x_1>0,x_n>0\},&&\\  
0,\quad\forall~x\in\partial 
B_1\cap\{x|x_1<0,x_n>0\},&&\\ 
-b,\quad\forall~x\in\partial 
B_1\cap\{x|x_n<0\}.&& 
\end{array} 
\right. 
\] 
Define 
\begin{equation} 
\phi(x):=\frac{1-|x|^2}{n\omega_n }\int_{\partial B_1}\frac{\phi_b 
(y)}{|x-y|^n}~dS_y,\quad\forall~x\in B_1,  
\label{abc+} 
\end{equation} 
where $\omega_n$ denotes the volume of the unit ball of $\Bbb{R}^n$.\newline
We know that $\phi\in C^{\infty}(B_1)\cap C^0(\bar 
B_1\setminus\partial\Bbb{R}^n_{+})$ and, after fixing some small $b>0$, 
\begin{equation} 
\left\{ 
\begin{array}{lcl} 
\Delta \phi=0\quad\mbox{in}~B_1,&&\\ 
\phi (0)>0,\quad 
\|\phi\|_{L^{\infty}(B_1)}\le 1,&&\\ 
\limsup\limits_{B_1\ni x\to\bar x}\phi (x)\le\max 
\{\bar {x}_1,0\},\quad\forall~\bar x\in\partial 
B_1 .&& 
\end{array} 
\right. 
\label{onething+} 
\end{equation} 
{\bf Claim.}\ There exists a constant $\tilde C>0$, depending only on 
$n,b$, such that  
\begin{equation} 
\frac{\partial\phi}{\partial x_n}(x)\ge \tilde C>0,\quad\forall~ x\in 
B_1\cap\partial\Bbb{R}^n_{+}.  
\label{twothings+} 
\end{equation} 
Indeed, consider 
\[
\psi (x',x_n):=\phi (x',x_n)-\phi (x',-x_n),\quad \eta (x):=\frac
b2 x_n,\quad\forall~x=(x',x_n)\in B_1^+. 
\]
We have
\[
\psi\ge\eta\quad\mbox{on}~\partial' B^+_1\cup (B_1\cap\partial\Bbb {R}^n_+).
\]
And for any $x\in B^+_1$, 
\begin{eqnarray*}
\psi (x)&=&\frac{1-|x|^2}{n\omega_n }\int_{\partial B_1}\frac{\phi_b 
(y)-\phi_b (y',-y_n)}{|x-y|^n}~dS_y\\
&=&\int_{\partial' B_1^+}\Big (\phi_b 
(y)-\phi_b (y',-y_n)\Big )\Big (\frac
{1}{|x-y|^n}-\frac{1}{|x-(y',-y_n)|^n}\Big )~dS_y>0, 
\end{eqnarray*}
therefore 
\[
\liminf_{x\to\tilde x}(\psi-\eta)(x)\ge 0,\quad\forall~\tilde x\in\partial
B_1^+.
\]
By Maximum Principle, $\psi\ge\eta$ in $B_1^+$. Since $\psi-\eta=0$ on
$B_1\cap\partial\Bbb{R}^n_+$, we have 
\[
\frac {\partial \psi}{\partial x_n}\ge\frac {\partial \eta}{\partial
x_n}=\frac b2,\quad\mbox{on}~B_1\cap\partial\Bbb{R}^n_+.
\]
The Claim is proved.\newline
{\bf Proof of Lemma~\ref{lemA4+}.}\ We only need to prove the 
lemma with $a=1, p-q=e_1:=(1\cdots,1)$. Indeed, replacing $u$ by $\frac 1a 
u$, $c$ by $ca^{\frac{2}{n-2}}$, $p$ by $\frac 1a p$ and $q$ by $\frac 
1a q$, we can assume $a=1$. After a rotation, we can assume 
$p-q=\lambda e_1$ for some $\lambda>0$. Replacing $u(x)$ by 
$u(\frac{x}{\lambda})$, $c$ by $\frac{c}{\lambda}$, $p,q$ by 
$\frac{p}{\lambda}, \frac{q}{\lambda}$ respectively, we can also  
assume $p-q=e_1$. \newline 
Since $\lim\limits_{r\to 0}\frac{\bar\delta (r)}{r}=0$, there exists 
$0<\bar r<d$ such that  
\begin{equation} 
\frac{\bar \delta (r)}{r}\le \frac 12\phi (0),\quad \forall~0<r<\bar r,
\label{cde+} 
\end{equation} 
where $\phi$ is defined by (\ref{abc+}).\newline
For $0<r<\bar r$, we consider, for $0<s<r$, 
\[ 
\phi^r(x):= 1+cx_n+q\cdot x'+r\phi (\frac xr)-\frac{ 
s^{n-2}d}{|x|^{n-2}}- \sup\limits_{(0,r]}\bar\delta,\quad\forall~x\in 
\overline{B_r\setminus B_s}. 
\] 
By the equations of $u$ and $\phi$, we have  
\[
\Delta (u-\phi^r)\le 0\quad\mbox{in}~B_r^{+}\setminus B_s^{+}. 
\] 
By the last lines in (\ref{A4e1+}) and (\ref{onething+}), 
\begin{equation} 
\limsup\limits_{B_r^{+}\ni x\to\bar x}(u(x)-\phi^r(x))\ge 
0,\quad\forall~\bar x\in \partial B_r\cap\overline{\Bbb{R}^n_{+}}. 
\label{f2+} 
\end{equation} 
Indeed, if $\bar x_1\ge 0$, we have, using $p-q=e_1$,
\[
u(x)-\phi^r (x)\ge (1+p\cdot x'+cx_n-\bar\delta (|x|))-\phi^r (x)\ge
x_1-r\phi (\frac xr),
\]
from which we deduce (\ref{f2+}).\newline
If $\bar x_1<0$, estimate (\ref{f2+}) follows from 
\[
u(x)-\phi^r (x)\ge (1+q\cdot x'+cx_n-\bar\delta (|x|))-\phi^r (x)\ge
-r\phi (\frac xr).
\]
Since $\|\phi\|_{L^{\infty}(B_1)}\le 1$, we have  
\[
\phi^r(x)<1+cx_n+q\cdot 
x'-\sup\limits_{(0,r]}\bar\delta,\quad\forall~x\in\partial 
B_s\cap\overline{\Bbb{R}^n_{+}}.  
\] 
Thus, by the last line in (\ref{A4e1+}), 
\begin{equation} 
u-\phi^r\ge 0\quad\mbox{on}~\partial B_s\cap\overline{\Bbb{R}^n_{+}}. 
\label{f3+} 
\end{equation} 
{\bf Claim.}\ There exists $\tilde r\in (0,\bar r]$, s.t., $\forall~ 
0<s<r<\tilde r$, 
\begin{equation} 
\inf\limits_{B_r^{+}\setminus B_s^{+}}(u-\phi^r)\ge 0. 
\label{claim'''} 
\end{equation}  
Suppose not, we have, by (\ref{f2+}), (\ref{f3+}), and the strong
maximum principle,   
\[ 
\inf\limits_{B_{r}^{+}\setminus 
B_s^{+}}(u-\phi^r)=(u-\phi^r)(\bar x)<0\quad\mbox{for some}~\bar x\in
(\partial\Bbb {R}^n_{+}\cap(B_{r}\setminus \overline{B_s})).   
\]  
At $\bar x$,  
\begin{equation} 
0\le \frac{\partial (u-\phi^r)}{\partial x_n}=cu^{\frac{n}{n-2}}(\bar
x)-c-\frac{\partial\phi}{\partial x_n}(\frac{\bar x}{r})\le c
u^{\frac{n}{n-2}}(\bar x)-c-\tilde C, 
\label{f6+} 
\end{equation} 
where $\tilde C$ is the constant in (\ref{twothings+}).\newline 
By the last line in (\ref{A4e1+}), we have, for some universal positive 
constant $C$,  
\[ 
u(\bar x)\ge 1-C|\bar x|\ge 1-Cr. 
\] 
On the other hand, 
\[ 
u(\bar x)\le \phi^r(\bar x)\le 1+Cr. 
\] 
We deduce from (\ref{f6+}), using the above two estimates,  
\[ 
0\le Cr-\tilde C, 
\] 
which is impossible if we choose $\tilde r<\min\{\frac{\tilde 
C}{C},\bar r\}$. (\ref{claim'''}) is established.\newline 
Sending $s\to 0$ in (\ref{claim'''}), we obtain 
\[ 
u(x)\ge 1+cx_n+q\cdot x'+r\phi (\frac 
xr)-\sup\limits_{(0,r]}\bar \delta,\quad\forall  
x\in B^{+}_r.  
\] 
Sending $x\to 0$, we have, by $(\ref{cde+})$, 
\[ 
\liminf\limits_{B_1^{+}\ni x\to 0} u(x)\ge 1+r\phi 
(0)-\sup\limits_{(0,r]}\bar \delta>1. 
\] 
Lemma~\ref{lemA4+} is established. 
\vskip 5pt 
\hfill $\Box$ 
\vskip 5pt 

\section{Appendix A}
\begin{lem} Let $a>0$ be a positive number and
 $\alpha$ be a real number.  Assume
 that  $h\in C^1[-4a,4a]$ satisfies, $\forall |\tau|<2a,~|s|\le
4a,~0<\lambda<a,~\lambda<|s-\tau|$,
\begin{equation}
\Big(\frac{\lambda}{|s-\tau|}\Big)^{\alpha}
h\Big(\tau+\frac{\lambda^2(s-\tau
)}{|s-\tau|^2}\Big)\le h(s).
\label{l1e1}
\end{equation}
Then
\[
|h^{\prime}(s)|\le \frac \alpha {2a}  h(s),\quad\forall\ |s|\le a.
\]
\label{al1}
\end{lem}  
{\bf Proof of Lemma~\ref{al1}.}\  By considering $h(as)$, we only need
to prove the lemma for $a=1$. If $\alpha=0$, it is easy to see that 
$h$ is identically equal to a constant on $[-1,1]$.
So we always assume that $\alpha\ne 0$.  We Only need to show that 
\begin{equation}
-h^{\prime}(s)\le  \frac \alpha 2  h(s),\quad\forall |s|<1,
\label{A1}
\end{equation}
since the estimate for $h^{\prime}(s)$ can be
obtained by applying  the above  $h(-s)$.

Now for $|\tau|<2$, let $h_{\tau}(s):=h(\tau+s)$, (\ref{l1e1}) is
equivalent to 
\[
\Big(\frac{\lambda}{|s-\tau|}\Big)^{\alpha}h_{\tau}\Big(\frac{\lambda^2(s-\tau
)}{|s-\tau|^2}\Big)\le h_{\tau}(s-\tau),\quad\forall |\tau|<2,~|s|\le
4,~0<\lambda<1,~\lambda<|s-\tau|,
\]
which implies, by setting $x=s-\tau$, that
$$
\Big(\frac{\lambda}{|x|}\Big)^{\alpha}h_{\tau}\Big(\frac{\lambda^2 x}
{|x|^2}\Big)\le h_{\tau}(x),\quad\forall |\tau|<2,
~0<\lambda<1,~\lambda<x<2.
$$
Let  $y=\frac{\lambda^2 x}{|x|^2}=\frac{\lambda^2}{x}$ in the above,
we have
\[
y^{\frac{\alpha}{2}}h_{\tau}(y)\le
x^{\frac{\alpha}{2}}h_{\tau}(x),\quad\forall\ 0<y<x<1.
\]
Thus 
\[
0\le\frac{d}{dx}\Big(x^{\frac{\alpha}{2}}h_{\tau}(x)\Big)=\frac{\alpha}{2}
x^{\frac{\alpha}{2}-1}h_{\tau}(x)+x^{\frac{\alpha}{2}}h_{\tau}^{\prime}(x),
\quad \forall \ 0<x<1,
\]
i.e
$$
\frac{\alpha}{2}h_{\tau}(x)+xh_{\tau}^{\prime}(x)\ge 0,\quad\forall 0<x<1.
$$
Let  $x\to 1$ in the above, we have 
\[
\frac{\alpha}{2}h_{\tau}(1)\ge -h_{\tau}^{\prime}(1),
\]
i.e.
\[
\frac{\alpha}{2}h(\tau+1)\ge -h^{\prime}(\tau+1),\quad\forall |\tau|<2.
\]
Estimate (\ref{A1}) follows from the above.

\vskip 5pt
\hfill $\Box$
\vskip 5pt

\begin{lem} Let $a>0$ be a constant and
let  $B_{8a}\subset \Bbb R^n$ be the ball of radius $8a$ and centered at the
origin, $n\ge 3$.
Assume that $u\in C^1(B_{8a})$ is a non-negative function
satisfying
$$
u_{x,\lambda}(y)\le u(y),\quad\forall \
x\in B_{4a}, \ y\in B_{8a},~0<\lambda<2a,~\lambda<|y-x|,
$$
where
$u_{x,\lambda}(y):=\Big(\frac{\lambda}{|y|}\Big)^{n-2}u\Big(x+\frac
{\lambda^2(y-x)}{|y-x|^2}\Big)$.
Then $\exists C(n)>0$, s.t.
\[
|\nabla u(x)|\le  \frac {n-2}{2a} u(x),\quad\forall |x|<a.
\]
\label{l2}
\end{lem}
{\bf Proof of Lemma~\ref{l2}.}\ For $\forall x\in B_a, e\in\Bbb R^n, |e|=1$,
let 
$h(s):=u(x+se)$. Then, by the hypothesis on $u$,
 $h$ satisfies the hypothesis of 
Lemma \ref{al1}.  Thus we have
$$
|h'(0)|\le \frac {n-2}{2a}h(0),
$$
i.e.
$$
|\nabla u(x)\cdot e|\le \frac {n-2}{2a} u(x).
$$
Lemma \ref{l2} follows from the above.

\vskip 5pt
\hfill $\Box$
\vskip 5pt

\section{Appendix B}
We first show that we may assume without loss of generality that the
$f$ in Theorem~\ref{yamabe}$^{\prime}$ is in addition homogeneous of
degree $1$. We achieve this by constructing the $\tilde f$ which is
homogeneous of degree $1$, $\tilde f^{-1} (1)=f^{-1} (1)$, and satisfies
the hypotheses of Theorem~\ref{yamabe}$^{\prime}$. \newline
By the cone structure of $\Gamma$,
the ray $\{s\lambda\ |\ s>0\}$ belongs to $\Gamma$ for
every $\lambda\in \Gamma$.
By the concavity of $f$,  we deduce from
(\ref{aa6}) that
\begin{equation}
\sum_{i=1}^nf_{\lambda_i}(\lambda) \lambda_i
>0, \qquad \forall\ \lambda\in \Gamma.
\label{aa6new}
\end{equation}
Since $f(0)=0$, $f$ satisfies (\ref{aa6}) and 
(\ref{aa6new}), and $f\in C^{4,\alpha}(\Gamma)$,  the equation
\begin{equation}
f(\varphi(\lambda)\lambda)=1, \qquad \lambda\in \Gamma
\label{ephi}
\end{equation}
defines, using the implicit function theorem, a positive function
$\varphi\in C^{4,\alpha}(\Gamma)$. It is easy to see from the
definition of $\varphi$ that
$\varphi(s\lambda)=s^{-1}\varphi(\lambda)$ for all $\lambda\in \Gamma$
and $0<s<\infty$. Set 
$$
\tilde f=\frac 1 \varphi, \qquad \mbox{on}\ \Gamma.
$$
By the homogeneity of $\varphi$, $\tilde f$   is homogeneous of
degree $1$. We will show that $\tilde f$ has the desired
properties. Clearly, $\tilde f$ is symmetric, (\ref{aa6}) is satisfied
and $\tilde f^{-1}(1)= f^{-1}(1)$.\newline 
To prove $\nabla\tilde f\in\Gamma_n$, applying
$\frac{\partial}{\partial \lambda_i}$ to (\ref{ephi}), we have
$$
0=f_{\mu_i}(\mu) \varphi(\lambda)+
\frac {\varphi_{\lambda_i}(\lambda)}{ \varphi(\lambda)}
\sum_{j=1}^n  f_{\mu_j}(\mu)\mu_j,
$$
where $\mu=\varphi(\lambda)\lambda$.
Since $f_{\mu_i}(\mu)>0$ and $\sum_{j=1}^n  f_{\mu_j}(\mu)\mu_j>0$, we have
$\varphi_{\lambda_i}(\lambda)<0$, i.e.,
$$
\tilde f_{\lambda_i}>0\qquad \mbox{on}\ \Gamma,\qquad \forall\ 1\le i\le n.
$$
Next we prove the concavity of $\tilde f$.   For $\lambda, \bar
\lambda\in \Gamma$, we have, by the concavity of $f$, that
\begin{eqnarray*}
&&f\Big(\frac{\varphi(\lambda)\varphi(\bar \lambda)}{t
\varphi(\bar \lambda)+(1-t)\varphi(\lambda)}
[t\lambda+(1-t)\bar\lambda]\Big)
\\
&=&
f\Big( \frac{t\varphi(\bar\lambda)}
            {t\varphi(\bar \lambda)+(1-t)\varphi(\lambda)}
        \varphi(\lambda)\lambda
+\frac{(1-t)\varphi(\lambda)}
            {t\varphi(\bar \lambda)+(1-t)\varphi(\lambda)}
        \varphi(\bar\lambda)\bar\lambda\Big)\\
&\ge &
\frac{t\varphi(\bar\lambda)}
            {t\varphi(\bar \lambda)+(1-t)\varphi(\lambda)}
f(\varphi(\lambda)\lambda)+
\frac{(1-t)\varphi(\lambda)}
            {t\varphi(\bar \lambda)+(1-t)\varphi(\lambda)}
f(\varphi(\bar\lambda)\bar\lambda)\\
&=&1=f\left(\varphi(t\lambda+(1-t)\bar\lambda)[t\lambda+(1-t)\bar\lambda]
\right).
\end{eqnarray*}
By (\ref{aa6new}), $f$ is strictly increasing along any ray in $\Gamma$
starting from the origin, therefore we deduce from the above that
$$
\frac{\varphi(\lambda)\varphi(\bar \lambda)}{t
\varphi(\bar \lambda)+(1-t)\varphi(\lambda)}
\ge \varphi(t\lambda+(1-t)\bar\lambda),
$$
i.e.,
$$
t\tilde f(\lambda)+(1-t)\tilde f(\bar \lambda)
\le \tilde f(t\lambda+(1-t)\bar \lambda).
$$
We have showed that $\tilde f$ is a concave function
in $\Gamma$.\newline
To check $\tilde f\in C^0(\overline \Gamma)$
and $\tilde f=0$ on $\partial \Gamma$, we only need to show that
$$
\lim_{\lambda\to \bar\lambda, \lambda\in \Gamma}
\tilde f(\lambda)=0\qquad\forall\ \bar\lambda\in \partial \Gamma.
$$
We show the above by contradiction argument.
Suppose the contrary, then for
some $\bar\lambda\in \partial \Gamma$ there exists a sequence
$\lambda^i\in \Gamma$, $\lambda^i\to \bar\lambda$, such that
$\lim_{i\to\infty}\tilde f(\lambda^i)>0$.  It follows that
$\varphi(\lambda^i)\to a$ for some $a\in[0,\infty)$.  By the continuity of
$f$ on $\overline \Gamma$,  we have
$1=f(\varphi(\lambda^i)\lambda^i)\to f(a\bar\lambda)$.
Since $f=0$ on $\partial \Gamma$, we have $a>0$ and
$\bar\lambda\in \Gamma$, a contradiction. We have proved that the
$\tilde f$ has the desired properties.
\vskip 3pt
\hfill $\Box$
\vskip 3pt
\begin{prop}
Let $V$ be an open symmetric convex subset of $\Bbb
R^n$ with $\partial V\neq\emptyset$. Assume that 
\begin{equation}
\nu(\lambda)\in \Gamma_n,\
 \qquad\forall\ \lambda\in \partial V,
\label{nu1new}
\end{equation}
and
\begin{equation}
\nu(\lambda)\cdot \lambda>0, \qquad \forall\ \lambda\in \partial V,
\label{nu2new}
\end{equation}
where $\nu(\lambda)$ denotes the  unit inner normal
of a supporting plane of $V$ at $\lambda$.
Then
$\Gamma(V)$ as defined in (\ref{cone}) is an open symmetric convex
cone with vertex at the origin. Moreover,
\begin{equation}
\Gamma_n
\subset \Gamma(V)\subset
\Gamma_1,
\label{g1}
\end{equation}
and
\begin{equation}
\Gamma(V)=\{s\lambda\ |\ \lambda\in \partial V, s>0\}.
\label{g2}
\end{equation}
\label{thmB-1}
\end{prop}
\begin{rem}
No regularity assumption on $\partial V$ is needed.
\end{rem}
To prove Proposition~\ref{thmB-1}, we need the following lemma.
\begin{lem}  Let $V$ be as in Proposition \ref{thmB-1}.  Then\newline
\noindent (i)\ If $\lambda\in V$, then
$\{s\lambda|~s\ge1\}\subset V$.\newline
\noindent (ii)\ $0\notin\bar V$.\newline
\noindent (iii)\ If $\lambda\in \partial V$, then
$\{s\lambda|~-\infty<s<1\}\cap\bar V=\emptyset$
and $\{s\lambda|~s>1\}\subset V$.
\label{t4c1}
\end{lem}
\noindent{\bf Proof of Lemma ~\ref{t4c1}.}\ If (i) does not hold, then
there exists some $\lambda\in V$ and $\bar s>1$ such that $\bar s
\lambda\in \partial V $. By the convexity of $V$, we have 
\[
(\lambda-\bar s\lambda)\cdot \nu (\bar
s\lambda)\ge 0.
\] 
From which, we deduce, by $\bar s>1$, that $\bar s\lambda\cdot\nu(\bar
s\lambda)\le 0$, contradicting (\ref{nu2new}). (i) is established.\newline
If $0\in \bar V$, by (\ref{nu2new}), $0\notin\partial V$. Hence $0\in V$. 
Since $V$ is open,  an open neighborhood of
$0$ belongs to $V$  and therefore,  by (i), $V=\Bbb R^n$,
contradicting the fact that $\partial V\neq \emptyset$.
(ii) is  established.\newline
Let $\lambda\in \partial V$.  For $-\infty<s<1$, we have, by
(\ref{nu2new}), that $\nu(\lambda)\cdot
(s\lambda-\lambda)=(s-1)\nu(\lambda)\cdot\lambda<0$. Since
$\nu(\lambda)$ is an inner normal, $s\lambda\notin \bar V$. Thus we
have proved the first statement in (iii). Now we prove the second
statement in (iii). Let $\lambda\in\partial V$, we know from the first
statement of (iii) that $\{s\lambda~|~s>1\}\cap\partial
V=\emptyset$. So either $\{s\lambda~|~s>1\}\subset V$ or
$\{s\lambda~|~s>1\}\cap V=\emptyset$. Noticing the first case is what we
want to prove, we can assume the second case. Then, in view of the
first statement of (iii), the line $\{s\lambda~|~s\in R\}$ has no
intersection with $V$. It follows, see theorem~11.2 in \cite{R}, that
there is a supporting plane of $V$ containing the line
$\{s\lambda~|~s\in R\}$, and therefore $\nu (\lambda)\cdot\lambda=0$,
where $\nu (\lambda)$ denotes the unit inner normal to the supporting
plane, contradicting (\ref{nu2new}).  (iii) is established.
\vskip 5pt
\hfill $\Box$
\vskip 5pt
\noindent{\bf Proof of
Proposition \ref{thmB-1}.}\  It is easy to see that 
$\Gamma(V)$ is an open symmetric convex cone with vertex at the
origin. Now we prove that $\Gamma(V)\subset \Gamma_1$.\newline
For any $ \lambda=(\lambda_1,\cdots,\lambda_n)\in\Gamma (V)$,
let
$$
\left\{ 
\begin{array}{lcl}
\lambda^1&=&\lambda=(\lambda_1,\cdots,\lambda_n),\\
\lambda^2&=&(\lambda_2, \cdots, \lambda_n, \lambda_1),\\ 
&\vdots&\\
\lambda^n&=&(\lambda_n, \lambda_1,\cdots, \lambda_{n-1}).
\end{array}
\right.
$$
Since $\Gamma (V)$ is symmetric, $\lambda^i\in \Gamma (V)$, $\forall~1\le i\le n$.
By the convexity of $\Gamma(V)$,
\[
\bar \lambda:=\frac{1}{n}\sum_{i=1}^n \lambda^i
=\frac{\sigma_1(\lambda)}{n}e\in\Gamma(V),
\]
where $e=(1,\cdots,1)$, $\sigma_1(\lambda)=\sum_{i=1}^n
\lambda_i$.\newline
Let
\[
\bar s:=\inf\{s>0|~s\bar\lambda\in V\}.
\]
By (ii) in Lemma \ref{t4c1}, $\bar s>0$ and $\bar s\bar \lambda\in\partial V$.
Let $\nu(\bar s\bar \lambda)$ be the unit inner normal of
a supporting plane of $V$ at $\bar s\bar\lambda$, we have, by (\ref{nu2new}),
\[
0<\nu(\bar s\bar\lambda)\cdot(\bar s
\bar\lambda)=\frac{\bar s}{n}\sigma_1(\nu(\bar s\bar \lambda))
\sigma_1(\lambda).
\]
By (\ref{nu1new}), $\sigma_1(\nu(\bar s\bar \lambda))>0$, thus
$\sigma_1(\lambda)>0$, i.e. $\Gamma (V)\subset \Gamma_1$.\newline
Next we prove $\Gamma_n\subset\Gamma$ by contradiction argument. 
Suppose $\exists~\mu\in\Gamma_n\setminus\Gamma(V)$. Take any
$\lambda\in\Gamma(V)\subset\Gamma_1$. Consider the $2$-dimensional plane
$\Bbb{P}$ generated by $\mu$ and $\lambda$.
We know that $\Gamma(V)\cap \Bbb P$ lies on
one side of the line $\partial \Gamma_1\cap \Bbb P$.
So $\{s\mu\ |\ s\in \Bbb R\}\cap \Gamma(V)=\emptyset$ and
therefore 
$\Gamma(V)\cap \Bbb P$ stays on one side of
$\{s\mu\ |\ s\in \Bbb R\}$ in $\Bbb P$, i.e.
$$
\bigg\{\tilde \mu\in \Bbb P\ |\
\tilde \mu\cdot\big[ \lambda-(\lambda\cdot \frac{\mu}{|\mu|})
 \frac{\mu}{|\mu|}
\big]<0\bigg\}\cap\Gamma(V)
=\emptyset.
$$
Fix some $\tilde \mu\in \Gamma_n\cap \Bbb P$ such that
$$
\tilde \mu\cdot\big[ \lambda-(\lambda\cdot \frac{\mu}{|\mu|})
 \frac{\mu}{|\mu|}\big]<0.
$$
Then the line $\ell:=\{s\tilde \mu\ |\ s\in \Bbb R\}$
has no intersection with $\overline V\cap \Bbb P$.
Now parallelly moving $\ell$ towards
$\overline V\cap\Bbb{P}$ and  a first touching
of the moving line and $\overline V\cap\Bbb{P}$ must occur.
Let $\bar \ell$ denote the first touching line and
let $\bar \lambda\in \bar \ell\cap (\overline V\cap\Bbb{P})$.
Clearly $\bar \lambda\in
\partial V$ and $\bar \ell \cap V=\emptyset$.
So there exists
a supporting plane of $V$ at $\bar \lambda$ which contains $\bar
{\ell}$.  Let $\nu(\bar \lambda)$ denote
the unit inner normal of the supporting plane, then,
$\nu(\bar \lambda)\cdot \tilde \mu=0$, a contradiction to 
$\tilde \mu\in\Gamma_n$ and $\nu (\bar\lambda)\in\Gamma_n$ by
(\ref{nu1new}). Thus $\Gamma_n\subset\Gamma(V)$. (\ref{g1}) is established.

Let
$$
\tilde{\Gamma}(V):=\{s\lambda|~\lambda\in\partial V, s>0\}.
$$
Next we show that $\Gamma(V)=\tilde{\Gamma}(V)$.
For $\lambda\in V$, consider the ray
$\{s\lambda|~s>0\}$. Since $0\notin\bar V$,  we know that
\[
\bar s:=\inf\{s|~s\lambda\in V\}>0.
\]
By the openness of $V$ and the definition of $\bar s$,
$\bar s\lambda\in\partial V$. So $\lambda \in 
\tilde \Gamma(V)$.  We have showed that 
$\Gamma(V)\subset \tilde{\Gamma}(V)$.
On the other hand, by (ii) and (iii) of Lemma  \ref{t4c1},
$\tilde{\Gamma}(V)\subset \Gamma(V)$.  We have established (\ref{g2}).
Proposition \ref{thmB-1} is established.
\vskip 3pt
\hfill $\Box$
\vskip 3pt
In the following, we deduce the equivalence of Theorem~\ref{yamabe} and
Theorem~\ref{yamabe}$^{\prime}$.\newline
{\bf Theorem~\ref{yamabe}$ $ $\Rightarrow$
Theorem~\ref{yamabe}$^{\prime}$.}\ Let $V:=\{\lambda\in\Gamma
|~f(\lambda)>1\}$. By (\ref{aa5}) and (\ref{aa6}), $\Gamma
(V)=\Gamma$. By the concavity and symmetry of $f$, $V$ is open, symmetric and
convex. Clearly $\partial V=\{\lambda\in\Gamma
|~f(\lambda)=1\}\neq\emptyset$ is $C^{4,\alpha}$ and $\nabla f$ is
inner normal to $\partial V$. Therefore $\nabla f\in\Gamma_n$ implies
(\ref{nu1}). In the above, we have proved the concavity of $f$ and
(\ref{aa6}) forces $\nabla f(\lambda)\cdot\lambda>0$. Restricted onto
$\partial V$, we have (\ref{nu2}). Hence
Theorem~\ref{yamabe}$^{\prime}$ follows from Theorem~\ref{yamabe}.\newline  
{\bf Theorem~\ref{yamabe}$^{\prime}$ $\Rightarrow$
Theorem~\ref{yamabe}.}\ We only need to construct a pair
$(f,\Gamma)$ satisfying all the assumptions 
in Theorem~\ref{yamabe}$^{\prime}$ and $\{f=1\}=\partial V$. Let
$\Gamma:=\Gamma (V)$ as defined in (\ref{cone}). By
Proposition~\ref{thmB-1}, (\ref{gamma0}) and (\ref{gamma1}) hold for
$\Gamma$. Let $f(s\lambda):=s$ for any $s\ge 0$ and any $\lambda\in\partial
V$. By (\ref{g2}), $\Gamma=\{s\lambda|~\lambda\in\partial V,~s>0\}$. So
$f$ is well defined,
symmetric and $C^{4,\alpha}$ on $\Gamma$. It is easy to see from the
definition that $f$ is homogeneous of degree $1$, therefore (\ref{aa6}) follows
directly. To prove $f$ is concave, taking any two points $a\lambda$
and $b\mu$ in $\Gamma$, where $\lambda,\mu\in\partial V$ and
$a,b>0$. For any $0\le t\le 1$, since $\lambda,\mu\in\partial V$ and
$\bar V$ is convex, we have
$\bar\lambda:=\frac{ta}{ta+(1-t)b}\lambda+\frac{(1-t)b}{ta+(1-t)b}\mu\in\bar 
V$. Recall the definition of $f$, we have
$\frac{\bar\lambda}{f(\bar\lambda)}\in \partial V$. However, by 
(i) and (iii)
in Lemma~\ref{t4c1}, we know $s\bar\lambda\in V$ for any $s>1$,
therefore $\frac{1}{f(\bar\lambda)}\le 1$, i.e., $f(\bar\lambda)\ge
1$. From which, we deduce that
\begin{eqnarray*}
&&f(ta\lambda+(1-t)b\mu)\\
&=&(ta+(1-t)b)
f\Big(\frac{ta}{ta+(1-t)b}\lambda+\frac{(1-t)b}{ta+(1-t)b}\mu\Big)\\
&=&(ta+(1-t)b)f(\bar\lambda)\ge ta+(1-t)b=tf(a\lambda)+(1-t)f(b\mu).
\end{eqnarray*} 
$f$ is concave.\newline
Now the only assumption left to check is that $f$ can be
continuously extended to $\partial\Gamma$ and varnishes on $\partial
\Gamma$. To see this, take any sequence $\{\lambda^i\}$ in $\Gamma$ with
$\lambda^i\to \bar\lambda\in\partial\Gamma$. We need to show
$\lim\limits_{i\to\infty}f(\lambda^i)=0$. Suppose not, there exists
a subsequence of $\{\lambda^i\}$, still denoted by $\{\lambda^i\}$ such
that $f(\lambda^i)\ge \delta$ for some constant $\delta>0$. By the
definition of $f$, $\frac{\lambda^i}{f(\lambda^i)}\in\partial
V$. On the other hand, $\lambda^i\to\bar\lambda$ and
$f(\lambda^i)\ge\delta>0$ implies $\{\frac{\lambda^i}{f(\lambda^i)}\}$
stays in a bounded set of $\Bbb {R}^n$. Hence
$\frac{\lambda^i}{f(\lambda^i)}\to \mu$ for some
$\mu\in\Bbb{R}^n$. Noticing $\partial V$ is closed,
$\mu\in\partial V$. By (\ref{g2}),
$\{s\mu~|~s>0\}\subset\Gamma$. Recall $0\notin\partial V$ and 
$\lambda^i\to\bar\lambda$,  we have $f(\lambda^i)$ is uniformly
bounded. W.l.o.g., we can assume $f(\lambda^i)\to c_0>0$. It follows
that 
\[
\partial\Gamma\ni\bar\lambda
\leftarrow\lambda^i=f(\lambda^i)\frac{\lambda^i}{f(\lambda^i)}\to
c_0\mu\in\Gamma, 
\]
a contradiction. Theorem~\ref{yamabe} follows from
Theorem~\ref{yamabe}$^{\prime}$. 
\vskip 3pt
\hfill $\Box$
\vskip 3pt
In the rest of this section, we address Remark~{\ref{remark1-1}}. We assume that
$\partial V\in C^{2,\alpha}$, but the principle curvatures of
$\partial V$ are positive. Let $P_1:=\partial \Gamma_1$. After a
rotation of the axis system, $\partial V$ can be represented as the 
graph of  a
$C^{2,\alpha}$ function $\bar\phi$ defined on $
P_1\equiv \Bbb {R}^{n-1}$ satisfying   
\begin{equation}
(\nabla ^2 \bar\phi)>0\quad\mbox {on}~\Bbb {R}^{n-1}.
\label{800}
\end{equation}
$\Gamma_n$ in the new axis system is still an open convex cone, denoted by
$\tilde \Gamma_n$. The assumptions (\ref{nu1new}) and (\ref{nu2new})
are translated into
\begin{equation}
(-\nabla\bar\phi (y'),1)\in\tilde\Gamma_n,~~(-\nabla\bar\phi
(y'),1)\cdot (y',\bar\phi (y'))>0,\qquad\forall~y'\in\Bbb{R}^{n-1}.
\label {condition}
\end{equation}
In the following, all the functions are
defined on $\Bbb {R}^{n-1}$ if not specified. For $R>0,~\epsilon>0$, consider
\[
\phi_R^{\epsilon}:=\rho_R\bar{\phi}^{\epsilon}+(1-\rho_R)\bar\phi,
\]
where $\bar {\phi}^{\epsilon}$ is the smooth mollifier of $\bar\phi$
and $\rho_R$ is a radially symmetric cut-off function having value $1$
in $B_R$ and $0$ outside $B_{2R}$. Let $V_{R}^{\epsilon}$ be the set
above the graph of $\phi_R^{\epsilon}$. For any $R>0$,
$\phi_R^{\epsilon}$ is identically equal to $\bar\phi$ outside
$B_{2R}$ and $\phi_r^{\epsilon}\to\bar\phi$ in $C_{loc}^{2,\alpha}$ as
$\epsilon \to 0$, so, for some small $\epsilon=\epsilon (R)>0$,
(\ref{nu1new}) and (\ref{nu2new}) hold for
$\phi_R^{\epsilon}$. Noticing $\partial V_{R}^{\epsilon}$ coincides
with $\partial V$ when $|y'|\ge 2R$, therefore we can assume, for the
same small $\epsilon$, $\Gamma (V_{R}^{\epsilon})=\Gamma
(V)$. Back to $\Gamma (V)$, $V$ and $V_R^{\epsilon}$ define $f$ and
$f_R^{\epsilon}$ as homogeneous functions of degree $1$ in $\Gamma
(V_R^{\epsilon})=\Gamma (V)$ taking value 1 on $\partial V$ and
$\partial V_R^{\epsilon}$ respectively. $f_R^{\epsilon}$ satisfies all
the assumptions of $f$ assumed in Theorem~\ref{yamabe}$^{\prime}$. Now
taking a sequence $R_i\to\infty$ and taking $\epsilon_i>0$ such that 
(\ref{800}) and (\ref{condition}) hold for $\phi_{R_i}^{\epsilon_i}$. Let
$f_{R_i}^{\epsilon_i}$ be the corresponding function on $\Gamma
(V)$. We know $f_{R_i}^{\epsilon_i}$ satisfies all the assumptions of
$f$ in Theorem~\ref{yamabe}$^{\prime}$ and is smooth in any compact
subset of $\Gamma$ and $f_{R_i}^{\epsilon_i}\to f$ in
$C_{loc}^{2,\alpha}(\Gamma)$. 

Consider the equation 
\begin{equation}
f_{R_i}^{\epsilon_i}(\lambda (A_{u^{\frac{4}{n-2}}g}))=1,\quad\lambda
(A_{u^{\frac{4}{n-2}}g} )\in\Gamma,\quad\mbox{on}~M^n.
\label{R}
\end{equation}
Applying
 Theorem~\ref{yamabe}$^{\prime}$  to  $(f_{R_i}^{\epsilon_i},\Gamma)$,
we have, for any solution $u_i$ of the equation (\ref{R}),
\begin{equation}
\|u_i\|_{C^{2,\alpha}(M^n,g)}+\|u_i^{-1}\|_{C^{2,\alpha}(M^n,g)}\le C 
\label {200}
\end{equation}
for some constant  $C$ is independent of $i$ ---   this is clear from the
proof of Theorem~\ref{yamabe}$^{\prime}$. 
This implies $\lambda (A_{u_i^{\frac{4}{n-2}}g})$ stays in a compact subset of
$\Gamma$ independent of $i$. Hence for $i$ large enough,
$f_{R_i}^{\epsilon_i}$ is $C^{4,\alpha}$ in this compact subset and we
have, by (\ref{200}) and Schauder theory, that
\[
\|u_i\|_{C^{4,\alpha}(M^n,g)}\le C_i,
\]  
where $C_i$ is some constant may depending on $i$.\newline
Following the degree arguments at the end of the proof of
Theorem~\ref{yamabe}$^{\prime}$ and replacing $O_t^{\ast}$ by 
\begin{eqnarray*}
O_t^i&:=&\{u\in C^{4,\alpha}(M^n,g)|
\|u\|_{C^{2,\alpha}(M^n,g)}+\|u^{-1}\|_{C^{2,\alpha}(M^n,g)}\le
2C,\\
&&\|u\|_{C^{4,\alpha}(M^n,g)}\le 2C_i,\quad\lambda
(A_{u^{\frac{4}{n-2}}g} )\in\Gamma \},
\end{eqnarray*}
we can find a solution $u_i$ of (\ref{R}). Since $u_i$ is uniformly
bounded in $C^{2,\alpha}(M^n,g)$, after passing to a subsequence,
$u_i$ converges in $C^{2}(M^n,g)$ to some function $u$ in
$C^{2,\alpha}(M^n,g)$. Sending $i\to\infty$ in (\ref{R}), we have
\[
f(\lambda (A_{u^{\frac{4}{n-2}}g}))=1,\quad\lambda
(A_{u^{\frac{4}{n-2}}g} )\in\Gamma,\quad\mbox{on}~M^n.
\]

\end{document}